\documentstyle{article}

\newcommand{\bigzerou}{%
\smash{\lower1.7ex\hbox{\bg 0}}}
\newtheorem{theorem}{Theorem}
\newtheorem{prop}{Proposition}
\newtheorem{defi}{Definition}
\newtheorem{cor}{Corollary}

\newtheorem{Rem}{Remark}
\newtheorem{lem}{Lemma}

\newcommand{\ba}{\begin{eqnarray}}
\newcommand{\ea}{\end{eqnarray}}
\newcommand{\no}{\nonumber}

\def\d{{\partial}}
\newcommand{\mapright}[1]{%
\smash{\mathop{%
\hbox to 1.0cm{\rightarrowfill}}\limits^{#1}}}
\newcommand{\mapleft}[1]{%
\smash{\mathop{%
\hbox to 1.3cm{\leftarrowfill}}\limits^{#1}}}

\begin{document}
\title{
\begin{flushright}
  \begin{minipage}[b]{5em}
    \normalsize
    ${}$      \\
  \end{minipage}
\end{flushright}
{\bf Gauss-Manin System and the Virtual Structure Constants}}
\author{Masao Jinzenji\\
\\
\it Division of Mathematics, Graduate School of Science \\
\it Hokkaido University \\
\it  Kita-ku, Sapporo, 060-0810, Japan\\
{\it e-mail address: jin@math.sci.hokudai.ac.jp}}
\maketitle
\begin{abstract}
In this paper, we discuss some applications of 
Givental's differential equations
to enumerative problems on rational curves in
projective hypersurfaces. 
Using this method, we prove some of the conjectures on the structure constants 
of quantum cohomology of projective hypersurfaces, proposed in our previous 
article.
Moreover, we clarify the correspondence between the virtual structure 
constants and Givental's differential equations when the projective 
hypersurface is Calabi-Yau or general type. 
\end{abstract}
\section{Introduction}
The main ingredient of this paper is the well-known ordinary differential 
equation: 
\begin{equation}
\biggl((\d_{x})^{N-1}-k\cdot e^{x}\cdot (k\d_{x}+k-1)(k\d_{x}+k-2)\cdots(k\d_{x}+1)\biggr)w(x)=0,
\label{fun}
\end{equation}
with arbitrary $N$ and $k$.
First, we derive the solutions of (\ref{fun}) that can be expressed as 
asymptotic expansion around $x=-\infty$.
To this aim, it is convenient to introduce the following rational 
function in $z$: 
\begin{equation}
\phi_{d}^{N,k}(z):=\frac{(kd)!}{(d!)^{N}}\prod_{j=1}^{kd}(1+\frac{k}{j}z)
\prod_{j=1}^{d}(1+\frac{1}{j}z)^{-N}.
\end{equation}
Then we introduce the generating function of $\phi_{d}^{N,k}(z)$:
\begin{equation}
\phi^{N,k}(x,z):=\sum_{d=0}^{\infty}\exp(dx)\frac{(kd)!}{(d!)^{N}}\prod_{j=1}^{kd\
}(1+\frac{k}{j}z)
\prod_{j=1}^{d}(1+\frac{1}{j}z)^{-N}.
\end{equation}
Next, we introduce the power series in $\exp(x)$,  
\begin{equation}
w_{j}^{N,k}(x):=(\d_{z})^{j}\phi^{N,k}(x,z)|_{x=0}.
\label{w}
\end{equation}
If we multiply $\phi^{N,k}(x,z)$ by $\exp(zx)$,  
\begin{equation}
\Phi^{N,k}(x,z):=\exp(zx)\phi^{N,k}(x,z),
\end{equation}
$\Phi^{N,k}(x,z)$ satisfies,
\begin{equation}
\biggl((\d_{x})^{N-1}-k\cdot e^{x}\cdot (k\d_{x}+k-1)(k\d_{x}+k-2)\cdots(k\d_{x}+1)\biggr)
\Phi^{N,k}(x,z)=z^{N-1}\cdot\exp(zx).
\end{equation}
Thus, we obtain the following $N-1$ solution of (\ref{fun}): 
\begin{equation}
u_{j}^{N,k}(x):=\frac{1}{j!}(\d_{z})^{j}
(\Phi^{N,k}(x,z))|_{z=0},\;\;\;(j=0,1,\cdots, N-2).
\end{equation}
When  $N-k\geq 2$, these solutions are the generating functions of a certain 
type of two-point correlation functions of topological sigma model on $
M_{N}^{k} $ : the degree
$k$ hypersurface in ${\bf P}^{N-1}$. In \cite{bert2}, Bertram and Kley
showed that 
all the rational correlation functions (inserted operators are restricted to 
K\"ahler sub-ring) are reconstructed from these solutions 
in the $N-k\geq 2$ case. Up to now, we also know how to 
modify $u_{j}^{N,k}(x)$ to 
construct the corresponding generating function when $N-k=1$. 

In this paper, we take a different path to reconstruct small quantum 
cohomology rings (K\"ahler sub-rings) from (\ref{fun}). 
Our idea is very simple. We just look at 
the differential equation instead of looking at the solution. We will first 
 show that the Gauss-Manin system associated with 
the quantum K\"ahler sub-ring 
of $M_{N}^{k}$ have the same informations as the ones of (\ref{fun}) if 
$N-k\geq 2$. Precisely speaking, if we assume 
the topological selection rule, we can 
determine all the structure constants of the  quantum K\"ahler sub-ring 
of $M_{N}^{k}$ from (\ref{fun}) via the Gauss-Manin system. Conversely, we can 
derive (\ref{fun}) by usual reduction of the 
 Gauss-Manin system associated with the quantum 
K\"ahler sub-ring \cite{giv}. 
We can extend our discussion to the $N-k\leq 0$ case.
In this case, direct relation between (\ref{fun}) and  
the Gauss-Manin system associated with  
quantum K\"ahler sub-ring of $M_{N}^{k}$ is lost. But we can still 
construct a kind of Gauss-Manin system which is directly connected to 
(\ref{fun}). In the $N=k$ case, this Gauss-Manin system is nothing but the 
B-model used in the mirror computation. 

On the other hand, we conjectured the recursive formulas that 
evaluate the structure constants of the quantum K\"ahler sub-ring 
of $M_{N}^{k}$ 
in terms of the ones of $M_{N+1}^{k}$ when $N-k\geq 2$ \cite{cj}, \cite{6}. 
These recursive formula is strong enough to determine all the structure 
constants of $M_{N}^{k}$ in this region. 
Then we find that the above reconstruction 
process via the Gauss-Manin system is useful enough 
to give a proof of the recursive formulas. 
The proof of them is one of the main results of this paper.
In \cite{cj} and \cite{gene}, 
we also conjectured that
 the virtual structure constants, that are obtained from iterated 
use of these recursive formulas into the 
$N-k\leq 0$ region, can be regarded as analogue of the B-model in the 
mirror computation. 
 We then constructed 
the generalized mirror transformation, that evaluate the structure constants 
of the quantum K\"ahler sub-ring of $M_{N}^{k}\;\;(N-k\leq 0)$ from 
the virtual 
structure constants, up to some lower degrees of rational curves
\cite{gene}, \cite{vir}. 

Another result of this paper is the assertion 
that the virtual structure constants
are nothing but the structure constants of the Gauss-Manin system 
associated with (\ref{fun}). Thus, we find a stronger evidence of the 
existence of the analogue of the mirror theorem for the general type 
projective hypersurface.  

This paper is organized as follows. In Section 2, we first introduce our  
notation for the quantum K\"ahler sub-ring of projective hypersurfaces. 
Next, we overview the conjectures proposed in our previous article, some 
of which are proved in this paper. In Section 3, we first introduce 
the Gauss-Manin system associated with the quantum K\"ahler sub-ring of 
$M_{N}^{k}\;\;(N-k\geq 2)$ and explain how all the structure constants 
of the quantum K\"ahler sub-ring are reconstructed from the Gauss-Manin 
system and (\ref{fun}).   Next, we prove the 
recursive formulas for the structure constants of the quantum K\"ahler 
sub-ring of $M_{N}^{k}\;\;(N-k\geq 2)$, introduced in the previous section.
Lastly, we extend our discussion in the $(N-k\geq 2)$ case to the cases 
$N-k=1$, $N-k=0$ and $N-k<0$ and show that the virtual structure constants 
introduced in Section 2 is nothing but the structure constants of the 
Gauss-Manin system associated with (\ref{fun}).  
\section{Quantum K\"ahler Sub-ring of Projective Hypersurfaces}
\subsection{Notation} 
In this section, we introduce the quantum K\"ahler sub-ring 
of the quantum cohomology ring of a degree $k$ hypersurface in
${\bf P}^{N-1}$.
Let $M_{N}^{k}$ be a hypersurface of degree $k$ in ${\bf P}^{N-1}$.
 We denote by $QH^{*}_{e}(M_{N}^{k})$ the 
sub-ring of the quantum cohomology ring $QH^{*}(M_{N}^{k})$
generated by ${\cal O}_{e}$ induced from the K\"ahler form $e$ 
(or, equivalently the intersection $H\cap M_{N}^{k}$ between a hyperplane
class $H$ of ${\bf P}^{N-1}$ and $M_{N}^{k}$).
Additive basis of $QH_{e}^{*}(M_{N}^{k})$ is given by 
${\cal O}_{e^{j}}\;\;(j=0,1,\cdots,N-2)$, which is induced from 
$e^{j}\in H^{j,j}(M_{N}^{k})$.
 The multiplication rule of $QH^{*}_{e}(M_{N}^{k})$ 
is determined by the Gromov-Witten invariant of genus $0$ 
$\langle {\cal O}_{e}{\cal O}_{{e}^{N-2-m}}
{\cal O}_{{e}^{m-1-(k-N)d}}\rangle_{d,M_{N}^{k}}$ and
it is given as follows:
\begin{eqnarray}
 L_{m}^{N,k,d} &:=&\frac{1}{k}\langle {\cal O}_{e}{\cal O}_{{e}^{N-2-m}}
{\cal O}_{{e}^{m-1-(k-N)d}}\rangle_{d},\no\\
\no\\
{\cal O}_{e}\cdot 1&=&{\cal O}_{e},\nonumber\\
{\cal O}_{e}\cdot{\cal O}_{{e}^{N-2-m}}&=&{\cal O}_{{e}^{N-1-m}}+
\sum_{d=1}^{\infty}L_{m}^{N,k,d}q^{d}{\cal O}_{{e}^{N-1-m+(k-N)d}},\no\\
q&:=&\exp(t),
\label{gm}
\end{eqnarray}
where the subscript $d$ counts the degree of the rational curves
measured by $e$. Therefore,  $q=\exp(t)$ is the degree counting 
parameter. 
\begin{defi}
We call $L_{n}^{N,k,d}$ the structure constant of weighted degree $d$.
\end{defi}
Since $M_{N}^{k}$ is a complex $(N-2)$ dimensional manifold, we see that
a structure constant $L_{m}^{N,k,d}$
is non-zero only if the following condition is satisfied:
\begin{eqnarray}
&& 1\leq N-2-m\leq N-2, 1\leq m-1+(N-k)d\leq N-2,\no\\
&\Longleftrightarrow &\mbox{max.}\{0,2-(N-k)d\}\leq m \leq 
\mbox{min.}\{N-3,N-1-(N-k)d\}.
\label{sel}
\end{eqnarray}
We can rewrite (\ref{sel}) into the form: 
\begin{eqnarray}
(N-k\geq 2) &\Longrightarrow& 0\leq m \leq (N-1)-(N-k)d\no\\
(N-k=1,d=1)&\Longrightarrow& 1\leq m \leq N-3\no\\
(N-k=1,d\geq2)&\Longrightarrow& 0\leq m \leq N-1-(N-k)d\no\\
(N-k\leq 0)&\Longrightarrow& 2+(k-N)d\leq m \leq N-3.
\label{flasel}
\end{eqnarray}
From (\ref{flasel}), we easily see that the number of the non-zero
structure constants $L_{m}^{N,k,d}$ is finite except for the case of $N=k$.
Moreover, if $N\geq 2k$, the non-zero structure constants come only from
the $d=1$ part and the non-vanishing $L_{m}^{N,k,1}$  
is determined by $k$ and  
independent of $N$. 
The $N\geq 2k$ region is studied by Beauville \cite{beauville}, 
and his result plays 
the role of an initial condition of our discussion later.
Explicitly, they are given by the formula :
\begin{equation}
\sum_{n=0}^{k-1}L_{n}^{N,k,1}w^{n}=k\prod_{j=1}^{k-1}(jw+(k-j)),
\label{one}
\end{equation}
and the other $L_{n}^{N,k,d}$'s all vanishes.
In the $N=k$ case, the multiplication rule of $QH^{*}_{e}(M_{k}^{k})$ is
given as follows:
\begin{eqnarray}
{\cal O}_{e}\cdot 1&=&{\cal O}_{e},\nonumber\\
{\cal O}_{e}\cdot{\cal O}_{{e}^{k-2-m}}&=&
(1+\sum_{d=1}^{\infty}q^{d}L_{m}^{k,k,d}){\cal O}_{{e}^{k-1-m}}
\;\;(m=2,3,\cdots,k-3),\no\\
{\cal O}_{e}\cdot{\cal O}_{{e}^{k-3}}  &=&{\cal O}_{e^{k-2}}.
\label{calabi}
\end{eqnarray}
Hence it is useful to  
introduce the generating function of the structure constants 
of the Calabi-Yau hypersurface $M_{k}^{k}$:
\begin{equation}
L_{m}^{k,k}(e^{t}):=1+\sum_{d=1}^{\infty}L_{m}^{k,k,d}e^{dt}\;\;(m=2,\cdots 
,k-3).
\end{equation}
\subsection{Overview of the Results for Fano and Calabi-Yau Hypersurfaces
and Introduction of the Virtual Structure Constants}
Let us summarize the conjectures proposed in 
 \cite{cj}, \cite{6}, some of which will be proved in this paper. 
In \cite{cj}, 
we conjectured that the structure constants $L_{m}^{N,k,d}$ 
of $QH_{e}^{*}(M_{N}^{k})$ for $(N-k\geq 2)$  can be obtained by 
applying the recursive formulas 
which describe
$L_{m}^{N,k,d}$ in terms of $L_{m'}^{N+1,k,d'}\;\;(d'\leq d)$,
with the initial 
conditions of $L_{m}^{N,k,1}$ given by (\ref{one})
and $L_{m}^{N,k,d}=0\;\;(d\geq 2)$ in the $N\geq 2k$ region.
Let us introduce the construction of the recursive formulas given in \cite{6}.
First, we introduce the polynomial $Poly_{d}$ in 
$x,y,z_{1},z_{2},\cdots,z_{d-1}$ defined by the formula: 
\begin{eqnarray}
&&Poly_{d}(x,y,z_{1},z_{2},\cdots,z_{d-1})\no\\
&&=\frac{1}{(2\pi\sqrt{-1})^{d-1}}
\int_{C_{1}}\frac{dt_{1}}{t_{1}}\cdots
\int_{C_{d-1}}\frac{dt_{d-1}}{t_{d-1}}\prod_{j=1}^{d-1}
\biggl(\frac{(d-j)x+jy}{d}+\sum_{i=1}^{j}\frac{d-j}{d-i}t_{i}+
\sum_{i=j+1}^{d-1}\frac{j}{i}t_{i}+\no\\
&&z_{j}(\frac{(d-j)x+jy}{d}+\sum_{i=1}^{j}\frac{d-j}{d-i}t_{i}+
\sum_{i=j+1}^{d-1}\frac{j}{i}t_{i})/
(\frac{(d-j)x+jy}{d}+\sum_{i=1}^{j}\frac{d-j}{d-i}t_{i}+
\sum_{i=j+1}^{d-1}\frac{j}{i}t_{i}-z_{j})\biggr)\no\\
&&=\frac{d}{(2\pi\sqrt{-1})^{d-1}}
\int_{D_{1}}{du_{1}}\cdots
\int_{D_{d-1}}{du_{d-1}}\prod_{j=1}^{d-1}\biggl((\frac{1}{2u_{j}-u_{j-1}-u_{j+1}})
\cdot(u_{j}+z_{j}\frac{u_{j}}{u_{j}-z_{j}})\biggr)\no\\
&&=\frac{d}{(2\pi\sqrt{-1})^{d-1}}
\int_{D_{1}}{du_{1}}\cdots
\int_{D_{d-1}}{du_{d-1}}\prod_{j=1}^{d-1}\biggl(
\frac{(u_{j})^{2}}{(2u_{j}-u_{j-1}-u_{j+1})(u_{j}-z_{j})}\biggr),
\label{trial2} 
\end{eqnarray}
where we denote $x$ (resp. $y$ ) by $u_{0}$ (resp. $u_{d}$) in the last two 
lines.
In (\ref{trial2}), the path $D_{i}$ 
goes around both poles $u_{i}=\frac{u_{i-1}+u_{i+1}}{2}, u_{i}=z_{i}$. 
Next, let us consider the monomial 
$x^{d_{i_{0}}}z_{i_{1}}^{d_{i_{1}}}\cdots z_{i_{l-1}}^{d_{i_{l-1}}}
y^{d_{i_{l}}}$
$(\sum_{j=0}^{l}d_{i_{j}}=d-1)$, that appear in $Poly_{d}$,
associated with the following ordered partition of a positive integer $d$ 
\cite{bert}:
\begin{equation}
0=i_{0}<i_{1}<i_{2}<\cdots<i_{l-1}<i_{l}=d.
\end{equation}
Then  we prepare some elements in (a free 
Abelian group) ${\bf Z}^{l}$, which are 
determined for each monomial 
$x^{d_{i_{0}}}z_{i_{1}}^{d_{i_{1}}}\cdots z_{i_{l-1}}^{d_{i_{l-1}}}
y^{d_{i_{l}}}$ ,
as follows:
\begin{eqnarray}
\alpha&:=&(l-d,l-d,\cdots,l-d),\no\\
\beta&:=&(0,i_{1}-1,i_{2}-2,\cdots,i_{l-1}-l+1),\no\\
\gamma&:=&(0,i_{1}(N-k),i_{2}(N-k),\cdots,i_{l-1}(N-k)),\no\\
\epsilon_{1}&:=&(1,0,0,0,\cdots,0),\no\\
\epsilon_{2}&:=&(1,1,0,0,\cdots,0),\no\\
\epsilon_{3}&:=&(1,1,1,0,\cdots,0),\no\\
&&\cdots\no\\
\epsilon_{l}&:=&(1,1,1,1,\cdots,1).
\end{eqnarray}
Now we define $\delta=(\delta_{1},\cdots,\delta_{l})\in{\bf Z}^{l}$ 
by the formula: 
\begin{equation}
\delta:=\alpha+\beta+\gamma+\sum_{j=1}^{l-1}(d_{i_{j}}-1)\epsilon_{j}+
d_{i_{l}}\epsilon_{l}.
\label{delta}
\end{equation}
With these set-up, we state the following theorem:
\begin{theorem}
The recursive formulas are given as follows:
\begin{equation}
L_{n}^{N,k,d}=\phi(Poly_{d}),
\label{mr}
\end{equation}
where $\phi$ is a ${\bf Q}$-linear map from the ${\bf Q}$-vector
space of the homogeneous polynomials of degree $d-1$ in $x,y,z_{1},\cdots
,z_{d-1}$ to the ${\bf Q}$-vector
space of the weighted homogeneous polynomials of degree $d$ in 
$L_{m}^{N+1,k,d'}$. And it is  defined on the basis by: 
\begin{equation}
\phi(x^{d_{0}}y^{d_{d}}z_{i_{1}}^{d_{i_{1}}}\cdots z_{i_{l-1}}^{d_{i_{l-1}}})
=\prod_{j=1}^{l}L_{n+\delta_{j}}^{N+1,k,i_{j}-i_{j-1}}.
\end{equation}
\end{theorem}
The proof will be given in the next section.
 
The structure constant $L_{m}^{k,k,d}$ for a Calabi-Yau hypersurface 
does not obey the recursive formulas (\ref{mr}).    
Instead, we introduce here the virtual structure constants 
$\tilde{L}_{m}^{N,k,d}$
as follows.
\begin{defi}
Let $\tilde{L}_{m}^{N,k,d}$ be the rational number obtained by applying
the recursion formulas (\ref{mr}) 
 for arbitrary $N$ and $k$ with the initial condition
 $L_{n}^{N,k,1}\;\;(N\geq 2k)$ and 
$L_{n}^{N,k,d}=0\;\;(d\geq 2,\;\;N\geq 2k)$. 
\end{defi}
\begin{Rem}
In the $N-k\geq2$ region, $\tilde{L}_{m}^{N,k,d}=L_{m}^{N,k,d}$.
\end{Rem}
\begin{Rem}
$\tilde{L}_{m}^{N,k,d}$ is non-zero if $0\leq m \leq N-1-(N-k)d$, and we have 
infinite number of $\tilde{L}_{m}^{N,k,d}$'s when $N-k\leq 0 $.
\end{Rem}
\begin{defi}
We call $\tilde{L}_{n}^{N,k,d}$ the virtual structure constant of 
weighted degree $d$.
\end{defi}
We define the generating function of the virtual structure constants of 
the Calabi-Yau hypersurface $M_{k}^{k}$
as follows:
\begin{eqnarray}
\tilde{L}^{k,k}_{n}(e^{x})&:=&1+
\sum_{d=1}^{\infty}\tilde{L}^{k,k,d}_{n}e^{dx},\nonumber\\
&&(n=0,1,\cdots,k-1).
\label{str2}
\end{eqnarray}
In \cite{cj}, we observed that $\tilde{L}^{k,k}_{n}(e^{x})$ gives us the 
information of the B-model of the mirror manifold of $M_{k}^{k}$. 
In this paper, we prove the theorem:
\begin{theorem}
\begin{eqnarray}
&&\tilde{L}_{0}^{k,k}(e^{x})=w_{0}^{k,k}(x)=
\sum_{d=0}^{\infty}\frac{(kd)!}{(d!)^{k}}e^{dx},\no\\
&&\tilde{L}_{1}^{k,k}(e^{x})=
\partial_{x}(x+\frac{w_{1}^{k,k}(x)}{w_{0}^{k,k}(x)})=\partial_{x}
\biggl(x+(\sum_{d=1}^{\infty}\frac{(kd)!}{(d!)^{k}}\cdot 
(\sum_{i=1}^{d}\sum_{m=1}^{k-1}\frac{m}{i(ki-m)})e^{dx})/
(\sum_{d=0}^{\infty}\frac{(kd)!}{(d!)^{k}}e^{dx})\biggr),\no\\
\label{po}
\end{eqnarray}
where $w_{j}^{k,k}(x)$ is the function introduced in (\ref{w}).
\end{theorem}
Of course, we can extend the theorem (\ref{po}) to the 
general $\tilde{L}_{n}^{k,k}(e^{x})$ if we compare the 
$\tilde{L}_{n}^{k,k}(e^{x})$ with the B-model three point functions
in \cite{gmp}.
In particular, this theorem asserts that we 
can obtain the mirror map $t=t(x)$ used in the mirror computation without 
assuming the mirror conjecture.
\begin{equation}
t(x)=x+\int_{-\infty}^{x}dx'({\tilde{L}^{k,k}_{1}(e^{x'})}-1)
=x+\sum_{d=1}^{\infty}\frac{\tilde{L}^{k,k,d}_{1}}{d}e^{dx}.
\end{equation}
With the above theorem, we can construct 
the mirror transformation that transforms 
the virtual structure constants of the  Calabi-Yau hypersurface into 
the real ones as follows:
\begin{equation}
L^{k,k}_{m}(e^{t})=\frac{\tilde{L}^{k,k}_{m}(e^{x(t)})}
{\tilde{L}^{k,k}_{1}(e^{x(t)})}.\quad(m=2,\cdots,k-3) 
\label{mt}
\end{equation}
After some combinatorial computation, we can rewrite (\ref{mt}) into 
the following form: 
\begin{equation}
L^{k,k,d}_{n}=\sum_{m=0}^{d-1}\mbox{Res}_{z=0}(z^{-m-1}
\exp(-d\sum_{j=1}^{\infty}\frac{\tilde{L}_{1}^{k,k,j}}{j}z^{j}))
\cdot(\tilde{L}^{k,k,d-m}_{n}-\tilde{L}^{k,k,d-m}_{1}).
\label{schur}
\end{equation}
In {\cite{gene}}, we argued that this formula must have deep connection
with toric compactification of the moduli space of rational curves 
in ${\bf P}^{N-1}$. With this idea, we speculated that we can generalize 
the formula (\ref{schur}) to the $N-k<0$ case. In \cite{gene} and \cite{vir},
we gave some numerical evidence of this generalization up to some lower
degree of rational curves.
\section{Gauss-Manin System}
Let us first introduce the Gauss-Manin system associated 
with the quantum K\"ahler sub-ring of $M_{N}^{k}$:
\begin{defi}
We call the following rank 1 ODE for vector valued function 
$\psi_{m}(t),\;\;
(m=0,1,\cdots,N-2)$:
\begin{eqnarray}
\partial_{t}\psi_{N-2-m}(t)&=&\psi_{N-1-m}(t)+\sum_{d=1}^{\infty}
\exp(dt)\cdot L_{m}^{N,k,d}\cdot\psi_{N-1-m-(N-k)d}(t),\no\\
\partial_{t}\psi_{N-2}(t)&=&\sum_{d=1}^{\infty}
\exp(dt)\cdot L_{0}^{N,k,d}\cdot\psi_{N-1-(N-k)d}(t),
\label{gm1}
\end{eqnarray}
the Gauss-Manin system associated with the quantum K\"ahler sub-ring of 
$M_{N}^{k}$.
\end{defi}
This definition can be applied to any $M_{N}^{k}$.
\subsection{Fano case ($N-k\geq 2$)}
In this case, we already have the celebrated theorem of Givental \cite{giv}:
\begin{theorem}{\bf (Givental)}\\
If $N-k\geq 2$, (\ref{gm1}) can be reduced to the rank $N-1$ ODE for 
$\psi_{0}(t)$:
\begin{equation}
\biggl((\partial_{t})^{N-1}-k\cdot e^{t}\cdot(k\partial_{t}+k-1)\cdots
(k\partial_{t}+2)\cdot(k\partial_{t}+1)\biggr)\psi_{0}(t)=0.
\label{gh} 
\end{equation}
\end{theorem}
Conversely, we can determine all the structure constants of the quantum 
K\"ahler sub-ring explicitly using (\ref{gm1}) and (\ref{gh}) as the starting 
point.
\begin{cor}
The structure constants $L_{n}^{N,k,d}$ are fully reconstructed 
from (\ref{gh}).
In particular, we have,
\begin{equation}
\sum_{n=0}^{k-1}\tilde{L}_{n}^{N,k,1}w^{n}=k\cdot\prod_{j=1}^{k-1}(jw+(k-j)).
\end{equation} 
\end{cor}
{\it proof)}\\
Using some algebra, we can represent $\psi_{N-1-m}(t)$ in terms of 
$\psi_{0}(t)$ as the form:
\begin{equation}
\psi_{N-1-m}(t)=(\d_{t})^{N-1-m}\psi_{0}(t)-
\sum_{d=1}^{\infty}\exp(dt)\sum_{j=0}^{N-1-m-(N-k)d}\gamma^{N,k,d}_{m,j}
(\d_{t})^{N-1-m-(N-k)d-j}\psi_{0}(t),
\label{ano}
\end{equation}
when $1\leq m \leq N-1$. Moreover, we can obtain the ODE for $\psi_{0}(t)$
by introducing $\psi_{N-1}(t)$ formally, which satisfies 
\begin{equation}
\d_{t}\psi_{N-2}(t)=\psi_{N-1}(t)+\sum_{d=1}^{\infty}\exp(dt)\cdot 
L_{0}^{N,k,d}\cdot\psi_{N-1-(N-k)d}(t).
\label{formal}
\end{equation}
If we represent $\psi_{N-1}(t)$ as the form of (\ref{ano}), the ODE is just 
given by the equation:
\begin{equation}
\psi_{N-1}(t)=0.
\label{0}
\end{equation}
Substitution of (\ref{ano}) into (\ref{gm1}) leads us to the recursive formula
for $\gamma^{N,k,d}_{m,j}$, 
\begin{equation}
\gamma_{m,0}^{N,k,d}-\gamma_{m+1,0}^{N,k,d}=L_{m}^{N,k,d}
-\sum_{f+g=d}L_{m}^{N,k,f}\gamma_{m+(N-k)f,0}^{N,k,g},
\label{rec1}
\end{equation}
\begin{equation}
\gamma_{m,j}^{N,k,d}-\gamma_{m+1,j}^{N,k,d}
=d\cdot \gamma_{m+1,j-1}^{N,k,d}-\sum_{f+g=d}L_{m}^{N,k,f}\gamma_{m+(N-k)f,j}^{N,k,g}.
\label{rec2}
\end{equation}
Here, we introduce the generating function,
\begin{equation}
\gamma_{m}^{N,k,d}(w):=\sum_{j=0}^{N-1-(N-k)d}\gamma_{m,j}^{N,k,d}w^{j}.
\end{equation}
Then, the recursive formulas (\ref{rec1}) and (\ref{rec2}) are reduced to one 
recursive formula for $\gamma_{m}^{N,k,d}(w)$,
\begin{equation}
\gamma_{m}^{N,k,d}(w)=(1+dw)\gamma_{m+1}^{N,k,d}(w)+L_{m}^{N,k,d}-
\sum_{f+g=d}L_{m}^{N,k,f}\gamma_{m+(N-k)f}^{N,k,g}(w).
\label{crec}
\end{equation}
Multiplying (\ref{crec}) by $(1+dw)^{m}$ makes the recursive formula 
more tractable. 
\begin{equation}
(1+dw)^{m}\gamma_{m}^{N,k,d}(w)-(1+dw)^{m+1}\gamma_{m+1}^{N,k,d}(w)
=(1+dw)^{m} L_{m}^{N,k,d}-
\sum_{f+g=d}(1+dw)^{m}L_{m}^{N,k,f}\gamma_{m+(N-k)f}^{N,k,g}(w).
\label{crec2}
\end{equation}
We can easily solve (\ref{crec2}) inductively. The answer is given by the 
formula:
\begin{eqnarray}
&&\gamma_{m}^{N,k,d}(w)=\no\\
&&\sum_{l=1}^{d}(-1)^{l-1}\sum_{(d_{1},\cdots,d_{l})\in OP_{d}}
\sum_{j_{l}=m}^{N-1-(N-k)d}\cdots\sum_{j_{2}=m}^{j_{3}}\sum_{j_{1}=m}^{j_{2}}
\prod_{i=1}^{l}\biggl((1+(\sum_{n=i}^{l}d_{n})w)^{j_{i}-j_{i-1}}
\cdot L_{j_{i}+(N-k)(\sum_{n=1}^{i-1}d_{n})}^{N,k,d_{i}}\biggr),\no\\
\label{A}
\end{eqnarray}
where we formally identify $j_{0}$ with $m$ and denote by $OP_{d}$ the set of 
the ordered partitions of $d$; $\{(d_{1},d_{2},\cdots,d_{l})|\;\; 
d_{j}\geq1,\;\;\sum_{j=1}^{l}d_{j}=d\}$.
At this point, we look back at Theorem 3. It merely says that
\begin{eqnarray}
\gamma_{0}^{N,k,1}(w)&=&k\cdot\prod_{j=1}^{k-1}(k+jw),\no\\
\gamma_{0}^{N,k,d}(w)&=&0,\;\;\;(d\geq 2).
\label{const}
\end{eqnarray}
Hence we obtain from (\ref{A}),
\begin{equation}
\sum_{m=0}^{k-1}L_{m}^{N,k,1}(1+w)^{m}=k\cdot\prod_{j=1}^{k-1}(k+jw),
\label{lin}
\end{equation}
and,
\begin{eqnarray}
&&\sum_{m=0}^{N-1-(N-k)d}L_{m}^{N,k,d}(1+dw)^{m}=\no\\
&&\sum_{l=2}^{d}(-1)^{l}\sum_{(d_{1},\cdots,d_{l})\in OP_{d}}
\sum_{j_{l}=0}^{N-1-(N-k)d}\cdots\sum_{j_{2}=0}^{j_{3}}\sum_{j_{1}=0}^{j_{2}}
\prod_{i=1}^{l}\biggl((1+(\sum_{n=i}^{l}d_{n})w)^{j_{i}-j_{i-1}}
\cdot L_{j_{i}+(N-k)(\sum_{n=1}^{i-1}d_{n})}^{N,k,d_{i}}\biggr).\no\\
\label{conic}
\end{eqnarray}
Substitution of $\frac{z-1}{d}$ into $w$ leads us to the formulas:
\begin{equation}
\sum_{m=0}^{k-1}L_{m}^{N,k,1}z^{m}=k\cdot\prod_{j=1}^{k-1}((k-j)+jz),
\end{equation}
\begin{eqnarray}
&&\sum_{m=0}^{N-1-(N-k)d}L_{m}^{N,k,d}z^{m}=\no\\
&&\sum_{l=2}^{d}(-1)^{l}\sum_{(d_{1},\cdots,d_{l})\in OP_{d}}
\sum_{j_{l}=0}^{N-1-(N-k)d}\cdots\sum_{j_{2}=0}^{j_{3}}\sum_{j_{1}=0}^{j_{2}}
\prod_{i=1}^{l}\biggl((1+
(\sum_{n=i}^{l}d_{n})(\frac{z-1}{d}))^{j_{i}-j_{i-1}}
\cdot L_{j_{i}+(N-k)(\sum_{n=1}^{i-1}d_{n})}^{N,k,d_{i}}\biggr).\no\\
\label{cinic}
\end{eqnarray}
It is obvious that we can completely determine 
all the $L_{n}^{N,k,d}$'s 
by induction of $d$, because the r.h.s. includes only the $L_{m}^{N,k,d'}$'s 
with $d'<d$.   Q.E.D.

{\bf Example}\\
$M_{7}^{5}$ model\\
The corresponding Gauss-Manin system is given by,
\begin{eqnarray} 
\d_{t}\psi_{0}(t)&=&\psi_{1}(t),\no\\
\d_{t}\psi_{1}(t)&=&\psi_{2}(t)+a\cdot e^{t}\cdot \psi_{0}(t), \no\\
\d_{t}\psi_{2}(t)&=&\psi_{3}(t)+b\cdot e^{t}\cdot \psi_{1}(t), \no\\
\d_{t}\psi_{3}(t)&=&\psi_{4}(t)+c\cdot e^{t}\cdot \psi_{2}(t)+d\cdot e^{2t}
\cdot \psi_{0}(t), \no\\
\d_{t}\psi_{4}(t)&=&\psi_{5}(t)+b\cdot e^{t}\cdot \psi_{3}(t)+g\cdot e^{2t}
\cdot \psi_{1}(t), \no\\
\d_{t}\psi_{5}(t)&=&\;\;\;\;\;\;\;\;\;\;\;\;a\cdot e^{t}\cdot \psi_{4}(t)+
d\cdot e^{2t}\cdot \psi_{2}(t)+f\cdot e^{3t}\cdot \psi_{0}(t),
\label{ex1}
\end{eqnarray}
where we used a trivial equality $L_{m}^{N,k,d}=L_{N-1-(N-k)d-m}^{N,k,d}$.
We reduce (\ref{ex1}) into an ordinary equation for $\psi_{0}(t)$ and obtain,
\begin{eqnarray}
&&\biggl((\d_{t})^{6}-e^{t}\cdot((2a+2b+c)(\d_{t})^{4}+(4a+4b+2c)(\d_{t})^{3}+
(6a+3b+c)(\d_{t})^{2}+(4a+b)(\d_{t})+a)\no\\
&&+e^{2t}(a^2+b^2+2ab+2ac-2d-g)(\d_{t})^{2}+e^{2t}(2a^2+2b^2+4ab+4ac-4d-2g)
(\d_{t})\no\\
&&+e^{2t}(a^2+2ab+4ac-4d)+e^{3t}(2ad-a^{2}c-f)\biggr)\psi_{0}(t)=0.
\label{clear}
\end{eqnarray} 
Then Theorem 1 asserts that the equation (\ref{clear}) equals the equation:
\begin{equation}
\biggl((\d_{t})^{6}-e^{t}\cdot(3125(\d_{t})^{4}+6250(\d_{t})^{3}+
4375(\d_{t})^{2}+1250(\d_{t})+120)\biggr)\psi_{0}(t)=0.
\label{com}
\end{equation}
By comparing (\ref{clear}) with (\ref{com}), we obtain,
\begin{equation}
a=120,\;b=770,\;c=1345,\;d=211200,\;g=692500,\;f=31320000,
\end{equation}
which agree  with our previous results in \cite{num}.

This corollary enables us to prove Theorem 1. As the first step, we prove 
the following theorem:
\begin{theorem}
Let $\varphi$ be the recursive formula in (\ref{mr}) considered as a 
homomorphism from the polynomial ring of $L_{n}^{N,k,d}$ to 
the one of $L_{n}^{N+1,k,d}$. Then we have, 
\begin{equation}
\varphi(\gamma_{0}^{N,k,d}(w))=\biggl(\prod_{j=1}^{d-1}(1+jw)
\biggr)\cdot\gamma_{0}^{N+1,k,d}(w).
\end{equation}
\end{theorem}
{\it proof)}
From now on, we use another notation 
$0=i_{0}<i_{1}<i_{2}<\cdots<i_{l-1}<i_{l}=d$ of ordered partition. 
Correspondence to the previous notation $(d_{1},\cdots,d_{l})\in OP_{d}$
is given by $d_{j}=i_{j}-i_{j-1}$.
We denote by 
$f_{(i_{1}-i_{0},\cdots,i_{l}-i_{l-1})}^{d_{i_{1}}d_{i_{2}}\cdots d_{i_{l-1}}}
(x,y)$
the coefficient polynomial 
of $z_{i_{1}}^{d_{i_{1}}}z_{i_{2}}^{d_{i_{2}}}\cdots 
z_{i_{l-1}}^{d_{i_{l-1}}}$ in the generating polynomial $Poly_{d}$.
Using (\ref{trial2}), it is explicitly given as follows:
\begin{eqnarray}
&&f_{(i_{1}-i_{0},\cdots,i_{l}-i_{l-1})}^{d_{i_{1}}d_{i_{2}}\cdots d_{i_{l-1}}}
(x,y)=\sum_{j=0}^{d-1-\sum_{m=1}^{l-1}d_{i_{m}}}
a_{d_{i_{1}}d_{i_{2}}\cdots d_{i_{l-1}}j}(i_{1}-i_{0},\cdots,i_{l}-i_{l-1})
x^{j}y^{d-1-\sum_{m=1}^{l-1}d_{i_{m}}-j}\no\\
&&:=\frac{d}{(2\pi\sqrt{-1})^{d-1}}\int_{D_{1}}du_{1}\cdots
\int_{D_{d-1}}du_{d-1}\prod_{m=1}^{d-1}\frac{u_{m}}{(2u_{m}-u_{m-1}-u_{m+1})}
\cdot\prod_{n=1}^{l-1}\frac{1}
{u_{i_{n}}^{d_{i_{n}}}}\no\\
&&=\frac{d}{(2\pi\sqrt{-1})^{l-1}}\int_{D_{i_{1}}}du_{i_{1}}\cdots
\int_{D_{i_{l-1}}}du_{i_{l-1}}\prod_{j=1}^{l}
f_{(i_{j}-i_{j-1})}(u_{i_{j-1}},u_{i_{j}})\times\no\\
&&\times\prod_{j=1}^{l-1}\frac{1}{((i_{j+1}-i_{j-1})u_{i_{j}}-
(i_{j}-i_{j-1})u_{i_{j+1}}-(i_{j+1}-i_{j})u_{i_{j-1}})\cdot 
u_{i_{j}}^{d_{i_{j}}-1}}\prod_{j=1}^{l-2}(i_{j+1}-i_{j}).
\label{exr}
\end{eqnarray}
where we have introduced the following polynomial:
\begin{equation}
\sum_{j=0}^{d-1}a_{j}(d)x^{j}y^{d-1-j}:=
\prod_{j=1}^{d-1}(\frac{jx+(d-j)y}{d})=
f_{(d)}(x,y).
\label{block}
\end{equation}
The path $D_{j}$ in the second line of (\ref{exr}) 
goes around $u_{j}=\frac{u_{j-1}+u_{j+1}}{2}, u_{j}=0$ 
if $j\in \{i_{1},\cdots,i_{l-1}\}$ and $u_{j}=\frac{u_{j-1}+u_{j+1}}{2}$
otherwise. The last equality in (\ref{exr}) is obtained from integrating out 
the variable $u_{j}\;\;(j\in\{1,2,\cdots,d-1\}\setminus 
\{i_{1},\cdots,i_{l-1}\})$.

With this definition, we can write down the form of the recursive formula
as follows, 
\begin{eqnarray}
&&L_{n}^{N,k,d}=\sum_{l=1}^{d}\sum_{0=i_{0}<\cdots<i_{l}=d}
\sum_{(d_{i_{1}},\cdots,d_{i_{l-1}})}\sum_{j=0}^{d-1-\sum_{m=1}^{l-1}d_{i_{m}}}
a_{d_{i_{1}}d_{i_{2}}\cdots d_{i_{l-1}}j}(i_{1}-i_{0},\cdots,i_{l}-i_{l-1})
\times\no\\
&&\times\prod_{h=1}^{l}L^{N+1,k,i_{h}-i_{h-1}}_{n-\sum_{m=1}^{h-1}d_{i_{m}}-j
+i_{h-1}(N-k+1)}.
\label{xr}
\end{eqnarray}
On the other hand, we can rewrite $\gamma_{0}^{N,k,d}(w)$ using the notation
$0=i_{0}<\cdots<i_{l}=d$ into, 
\begin{eqnarray}
\sum_{l=1}^{d}(-1)^{l-1}\sum_{0=i_{0}<\cdots<i_{l}=d}
\sum_{j_{l}=0}^{N-1-(N-k)d}\sum_{j_{l-1}=0}^{j_{l}}\cdots\sum_{j_{1}=0}^{j_{2}}
\prod_{n=1}^{l}\biggl((1+(d-i_{n-1})w)^{j_{n}-j_{n-1}}
L_{j_{n}+(N-k)i_{n-1}}^{N,k,i_{n}-i_{n-1}}\biggr).\no\\
\label{exg}
\end{eqnarray}
Substituting (\ref{xr}) into (\ref{exg}), we obtain,
\begin{eqnarray}
&&\varphi(\gamma_{0}^{N,k,d}(w))=\no\\
&&\sum_{l=1}^{d}(-1)^{l-1}\sum_{l=1}^{d}\sum_{0=i_{0}<\cdots<i_{l}=d}
\sum_{j_{l}=0}^{N-1-(N-k)d}\sum_{j_{l-1}=0}^{j_{l}}\cdots\sum_{j_{1}=0}^{j_{2}}
\no\\
&&\prod_{n=1}^{l}\biggl((1+(d-i_{n-1})w)^{j_{n}-j_{n-1}}
\sum_{l_{n}=1}^{i_{n}-i_{n-1}}\sum_{i_{n-1}=i_{0}^{n}<\cdots
<i_{l_{n}}^{n}=i_{n}}
\sum_{(d_{i_{1}^{n}},\cdots,d_{i_{l_{n}-1}^{n}})}
\sum_{c_{n}=0}^{i_{n}-i_{n-1}-1-\sum_{m_{n}=1}^{l_{n}-1}d_{i_{m_{n}}^{n}}}
\times\no\\
&&\times a_{d_{i_{1}^{n}}\cdots d_{i_{l_{n}-1}^{n}}c_{n}}
(i_{1}^{n}-i_{0}^{n},\cdots,i_{l_{n}}^{n}-i_{l_{n}-1}^{n})
\prod_{h_{n}=1}^{l_{n}}
L^{N+1,k,i_{h_{n}}^{n}-i_{h_{n}-1}^{n}}_{j_{n}-c_{n}-i_{n-1}-
\sum_{m_{n}=1}^{h_{n}-1}d_{i_{m_{n}}^{n}}
+i_{h_{n}-1}^{n}(N-k+1)}\biggr).
\end{eqnarray}
In the above formula, we can see appearance of iterated ordered 
partition:
\begin{equation}
0=i_{0}=i_{0}^{1}<i_{1}^{1}<\cdots<i_{l_{1}}^{1}=i_{1}=i_{0}^{2}<\cdots
<i_{l_{2}}^{2}=i_{2}<\cdots<i_{l_{l}}^{l}=i_{l}=d.
\end{equation}
Then we pick up the terms whose iterated ordered partition is equal to 
the ordered  partition $0=i_{0}<i_{1}<\cdots<i_{l}=d$. The  result is 
conveniently 
written in terms of ordered partition $0=h_{0}<h_{1}<\cdots<h_{s}=l$, 
and  the statement of the theorem is reduced to the following equality: 
\begin{eqnarray}
&&\sum_{s=1}^{l}(-1)^{s-1}\sum_{0=h_{0}<\cdots<h_{s}=l}
\sum_{j_{s}=0}^{N-1-(N-k)d}\sum_{j_{s-1}=0}^{j_{s}}\cdots\sum_{j_{1}=0}^{j_{2}}
\sum_{{(d_{i_{1}},\cdots,d_{i_{l-1}})}\atop{(d_{i_{h_{j}}}=0)}}
\prod_{n=1}^{s}\biggl((1+(d-i_{h_{n-1}})w)^{j_{n}-j_{n-1}}\times\no\\
&&\times\sum_{c_{n}=0}^{i_{h_{n}}-i_{h_{n-1}}-1-\sum_{j=h_{n-1}}^{h_{n}}d_{i_{j}}}
 a_{d_{i_{h_{n-1}+1}}\cdots d_{i_{h_{n}-1}}c_{n}}
(i_{h_{n-1}+1}-i_{h_{n-1}},\cdots,i_{h_{n}}-i_{h_{n}-1})\times
\no\\
&&\times\prod_{a=h_{n-1}+1}^{h_{n}}
L^{N+1,k,i_{a}-i_{a-1}}_{j_{n}-c_{n}-i_{h_{n-1}}-\sum_{j=h_{n-1}+1}^{a-1}
d_{i_{j}}+i_{a-1}(N-k+1)}\biggr)\no\\
&&=\biggl(\prod_{j=1}^{d-1}(1+jw)\biggr)
\cdot\biggl((-1)^{l-1}
\sum_{j_{l}=0}^{N-(N-k+1)d}\sum_{j_{l-1}=0}^{j_{l}}\cdots\sum_{j_{1}=0}^{j_{2}}
\prod_{n=1}^{l}(1+(d-i_{n-1})w)^{j_{n}-j_{n-1}}
L_{j_{n}+(N-k+1)i_{n-1}}^{N,k,i_{n}-i_{n-1}}\biggr).\no\\
\label{exe}
\end{eqnarray}
Next, we carefully look at the summand in the l. h. s. of (\ref{exe}) coming 
from the ordered partition $0=h_{0}<h_{1}<\cdots<h_{s}=l$:
\begin{eqnarray}
&&(-1)^{s-1}\sum_{{(d_{i_{1}},\cdots,d_{i_{l-1}})},\;{(d_{i_{h_{j}}}=0)}}
\sum_{c_{1}=0}^{i_{h_{1}}-i_{h_{0}}-1-\sum_{j=h_{0}}^{h_{1}}d_{i_{j}}}
\cdots
\sum_{c_{s}=0}^{i_{h_{s}}-i_{h_{s-1}}-1-\sum_{j=h_{s-1}}^{h_{s}}d_{i_{j}}}\no\\
&&\prod_{n=1}^{s} a_{d_{i_{h_{n-1}+1}}\cdots d_{i_{h_{n}-1}}c_{n}}
(i_{h_{n-1}+1}-i_{h_{n-1}},\cdots,i_{h_{n}}-i_{h_{n}-1})\times
\sum_{j_{s}=0}^{N-1-(N-k)d}\sum_{j_{s-1}}^{j_{s}}\cdots\sum_{j_{1}=0}^{j_{2}}
\no\\
&&\prod_{n=1}^{s}\biggl((1+(d-i_{h_{n-1}})w)^{j_{n}-j_{n-1}}
\prod_{a=h_{n-1}+1}^{h_{n}}
L^{N+1,k,i_{a}-i_{a-1}}_{j_{n}-c_{n}-i_{h_{n-1}}-\sum_{j=h_{n-1}+1}^{a-1}
d_{i_{j}}+i_{a-1}(N-k+1)}\biggr).\no\\
\label{addend}
\end{eqnarray}
Changing $j_{n}$ into $j_{n}'=j_{n}-c_{n}-i_{h_{n-1}}+
\sum_{j=1}^{h_{n-1}}d_{i_{j}}$, we can separate 
(\ref{addend}) into a bulk part:
\begin{eqnarray}
&&(-1)^{s-1}\sum_{{(d_{i_{1}},\cdots,d_{i_{l-1}})},\;{(d_{i_{h_{j}}}=0)}}
\sum_{c_{1}=0}^{i_{h_{1}}-i_{h_{0}}-1-\sum_{j=h_{0}}^{h_{1}}d_{i_{j}}}
\cdots
\sum_{c_{s}=0}^{i_{h_{s}}-i_{h_{s-1}}-1-\sum_{j=h_{s-1}}^{h_{s}}d_{i_{j}}}\no\\
&&\prod_{n=1}^{s} a_{d_{i_{h_{n-1}+1}}\cdots d_{i_{h_{n}-1}}c_{n}}
(i_{h_{n-1}+1}-i_{h_{n-1}},\cdots,i_{h_{n}}-i_{h_{n}-1})\times
\sum_{j_{s}=0}^{N-(N-k+1)d}\sum_{j_{s-1}=0}^{j_{s}}\cdots\sum_{j_{1}=0}^{j_{2}}
\no\\
&&\prod_{n=1}^{s}\biggl((1+(d-i_{h_{n-1}})w)^{j_{n}-j_{n-1}+c_{n}-c_{n-1}+
i_{h_{n-1}}-i_{h_{n-2}}-\sum_{j=h_{n-2}}^{h_{n-1}}d_{i_{j}}}
\prod_{a=h_{n-1}+1}^{h_{n}}
L^{N+1,k,i_{a}-i_{a-1}}_{j_{n}-\sum_{j=1}^{a-1}
d_{i_{j}}+i_{a-1}(N-k+1)}\biggr),\no\\
\label{bulk}
\end{eqnarray}
and boundary parts:
\begin{eqnarray}
&&(-1)^{s-1}\sum_{m=1}^{s-1}
\sum_{{(d_{i_{1}},\cdots,d_{i_{l-1}})},\;{(d_{i_{h_{j}}}=0)}}
\sum_{c_{1}=0}^{i_{h_{1}}-i_{h_{0}}-1-\sum_{j=h_{0}}^{h_{1}}d_{i_{j}}}
\cdots
\sum_{c_{s}=0}^{i_{h_{s}}-i_{h_{s-1}}-1-\sum_{j=h_{s-1}}^{h_{s}}d_{i_{j}}}
\no\\
&&\prod_{n=1}^{s} a_{d_{i_{h_{n-1}+1}}\cdots d_{i_{h_{n}-1}}c_{n}}
(i_{h_{n-1}+1}-i_{h_{n-1}},\cdots,i_{h_{n}}-i_{h_{n}-1})
\sum_{l_{m}=1}^{c_{m+1}-c_{m}+
i_{h_{m}}-i_{h_{m-1}}-\sum_{j=h_{m-1}}^{h_{m}}d_{i_{j}}}\no\\
&&\sum_{j_{s}=-c_{s}-i_{h_{s-1}}+l_{m}+
\sum_{j=1}^{h_{s-1}}d_{i_{j}}}^{N-1-(N-k)d-c_{s}-i_{h_{s-1}}+l_{m}+
\sum_{j=1}^{h_{s-1}}d_{i_{j}}}\cdots
\sum_{j_{m+2}=-c_{m+2}-i_{h_{m+1}}+l_{m}+
\sum_{j=1}^{h_{m+1}}d_{i_{j}}}^{j_{m+3}+c_{m+3}-c_{m+2}+i_{h_{m+2}}-
i_{h_{m+1}}-\sum_{j=h_{m+1}}^{h_{m+2}}d_{i_{j}}}\no\\
&&\sum_{j_{m}=0}^{j_{m+2}+
c_{m+2}-c_{m+1}+i_{h_{m+1}}-
i_{h_{m}}-\sum_{j=h_{m}}^{h_{m+1}}d_{i_{j}}}
\sum_{j_{m-1}=0}^{j_{m}}\sum_{j_{m-2}=0}^{j_{m-1}}
\cdots
\sum_{j_{1}=0}^{j_{2}}\no\\
&&\biggl(\prod_{n=1}^{m}(1+(d-i_{h_{n-1}})w)^{j_{n}-j_{n-1}}\biggr)\cdot
(1+(d-i_{h_{m+1}})w)^{j_{m+2}-j_{m}}\cdot
\biggl(\prod_{n=m+3}^{s}(1+(d-i_{h_{n-1}})w)^{j_{n}-j_{n-1}}\biggr)\no\\
&&\prod_{n=1}^{m}\biggl((1+(d-i_{h_{n-1}})w)^{c_{n}-c_{n-1}+
i_{h_{n-1}}-i_{h_{n-2}}-\sum_{j=h_{n-2}}^{h_{n-1}}d_{i_{j}}}
\prod_{a=h_{n-1}+1}^{h_{n}}
L^{N+1,k,i_{a}-i_{a-1}}_{j_{n}-\sum_{j=1}^{a-1}
d_{i_{j}}+i_{a-1}(N-k+1)}\biggr)\no\\
&&\biggl((1+(d-i_{h_{m}})w)^{-l_{m}
+c_{m+1}-c_{m}+
i_{h_{m}}-i_{h_{m-1}}-\sum_{j=h_{m-1}}^{h_{m}}d_{i_{j}}}
\prod_{a=h_{m}+1}^{h_{m+1}}
L^{N+1,k,i_{a}-i_{a-1}}_{j_{m}-l_{m}-\sum_{j=1}^{a-1}
d_{i_{j}}+i_{a-1}(N-k+1)}\biggr)\no\\
&&\prod_{n=m+2}^{s}\biggl((1+(d-i_{h_{n-1}})w)^{c_{n}-c_{n-1}+
i_{h_{n-1}}-i_{h_{n-2}}-\sum_{j=h_{n-2}}^{h_{n-1}}d_{i_{j}}}
\prod_{a=h_{n-1}+1}^{h_{n}}
L^{N+1,k,i_{a}-i_{a-1}}_{j_{n}-l_{m}-\sum_{j=1}^{a-1}
d_{i_{j}}+i_{a-1}(N-k+1)}\biggr).\no\\
\label{bound}
\end{eqnarray}
To derive (\ref{bound}), 
we further replace $j_{n}'\;(n>m)$ by $j_{n}''=j_{n}'+l_{m}$
and omit the dashes in $j_{n}', j_{n}''$ in the final form. 
Since we can easily see that the identity: 
\begin{equation}
\sum_{c_{1}=0}^{i_{1}-i_{0}-1}\cdots
\sum_{c_{l}=0}^{i_{l}-i_{l-1}-1}a_{c_{1}}(i_{1}-i_{0})\cdots
a_{c_{l}}(i_{l}-i_{l-1})\prod_{j=1}^{l}
\biggl((1+(d-i_{j-1})w)^{c_{j}-c_{j-1}+i_{j-1}-i_{j-2}}\biggr)
=\prod_{j=1}^{d-1}(1+jw),
\label{wow}
\end{equation}
holds true, the bulk part coming from the ordered partition 
$h_{j}=j\;\;(j=1,2,\cdots,l-1)$ is 
nothing but the r.h.s. of (\ref{exe}). 
Therefore, what remains to show is cancellation 
of the remaining terms. At this stage, 
we add some comments on boundary parts.
Looking at the first boundary part in (\ref{bound}) which 
corresponds to the operation to remove $m\;(m=1,2,\cdots,s-1)$ 
from the set $\{1,2,\cdots,s-1\}$,
we can further pick up the second boundary part which corresponds 
to remove $n\;(n=m+1,m+2,\cdots,s-1)$ from $\{m+1,m+2,\cdots,s-1\}$.
Explicitly, the first boundary part separated from the second boundary 
parts is given by the formula:
\begin{eqnarray}
&&(-1)^{s-1}
\sum_{{(d_{i_{1}},\cdots,d_{i_{l-1}})},\;{(d_{i_{h_{j}}}=0)}}
\sum_{c_{1}=0}^{i_{h_{1}}-i_{h_{0}}-1-\sum_{j=h_{0}}^{h_{1}}d_{i_{j}}}
\cdots
\sum_{c_{s}=0}^{i_{h_{s}}-i_{h_{s-1}}-1-\sum_{j=h_{s-1}}^{h_{s}}d_{i_{j}}}
\no\\
&&\prod_{n=1}^{s} a_{d_{i_{h_{n-1}+1}}\cdots d_{i_{h_{n}-1}}c_{n}}
(i_{h_{n-1}+1}-i_{h_{n-1}},\cdots,i_{h_{n}}-i_{h_{n}-1})
\sum_{l_{m}=1}^{c_{m+1}-c_{m}+
i_{h_{m}}-i_{h_{m-1}}-\sum_{j=h_{m-1}}^{h_{m}}d_{i_{j}}}\no\\
&&\sum_{j_{s-1}=0}^{N-(N-k+1)d}\sum_{j_{s-2}=0}^{j_{s-1}} 
\cdots\sum_{j_{1}=0}^{j_{2}}
\biggl(\prod_{n=1}^{m}(1+(d-i_{h_{n-1}})w)^{j_{n}-j_{n-1}}\biggr)\cdot
\biggl(\prod_{n=m+1}^{s-1}(1+(d-i_{h_{n}})w)^{j_{n}-j_{n-1}}\biggr)\no\\
&&\prod_{n=1}^{m}\biggl((1+(d-i_{h_{n-1}})w)^{c_{n}-c_{n-1}+
i_{h_{n-1}}-i_{h_{n-2}}-\sum_{j=h_{n-2}}^{h_{n-1}}d_{i_{j}}}
\prod_{a=h_{n-1}+1}^{h_{n}}
L^{N+1,k,i_{a}-i_{a-1}}_{j_{n}-\sum_{j=1}^{a-1}
d_{i_{j}}+i_{a-1}(N-k+1)}\biggr)\no\\
&&\biggl((1+(d-i_{h_{m}})w)^{-l_{m}
+c_{m+1}-c_{m}+
i_{h_{m}}-i_{h_{m-1}}-\sum_{j=h_{m-1}}^{h_{m}}d_{i_{j}}}
\prod_{a=h_{m}+1}^{h_{m+1}}
L^{N+1,k,i_{a}-i_{a-1}}_{j_{m}-l_{m}-\sum_{j=1}^{a-1}
d_{i_{j}}+i_{a-1}(N-k+1)}\biggr)\no\\
&&\prod_{n=m+2}^{s}\biggl((1+(d-i_{h_{n-1}})w)^{c_{n}-c_{n-1}+
i_{h_{n-1}}-i_{h_{n-2}}-\sum_{j=h_{n-2}}^{h_{n-1}}d_{i_{j}}}
\prod_{a=h_{n-1}+1}^{h_{n}}
L^{N+1,k,i_{a}-i_{a-1}}_{j_{n}-l_{m}-\sum_{j=1}^{a-1}
d_{i_{j}}+i_{a-1}(N-k+1)}\biggr).\no\\
\label{1stbound}
\end{eqnarray}
Continuing the same operation, we can observe that the summand 
of (\ref{addend}) produce $s-1 \choose t$ $t$-th boundary parts and 
that they have the same structure of summation on $j_{n}'s$ as the 
bulk part coming from the 
ordered partition:
\begin{eqnarray}
&&0=h_{p_{0}}<h_{p_{1}}<h_{p_{2}}<\cdots<h_{p_{s-t}}=l\no\\
&&(0=p_{0}<p_{1}<\cdots<p_{s-t}=s,\;\;1\leq t \leq s-1).
\label{boup}
\end{eqnarray}
We then separate the set $\{1,2,\cdots,s-1\}$ into disjoint union of two 
sets associated with (\ref{boup}).
\begin{eqnarray}
&&\{1,2,\cdots,s-1\}=\{p_{1},p_{2},\cdots,p_{s-t-1}\} \coprod
\{r_{1},r_{2},\cdots,r_{t}\}, \no\\
&& (r_{1}<r_{2}<\cdots<r_{t}).
\label{union}    
\end{eqnarray}
With these set-up, we can write down the $t$-th boundary part 
as the generalization of (\ref{1stbound}), 
\begin{eqnarray}
&&(-1)^{s-1}
\sum_{{(d_{i_{1}},\cdots,d_{i_{l-1}})},\;{(d_{i_{h_{j}}}=0)}}
\sum_{c_{1}=0}^{i_{h_{1}}-i_{h_{0}}-1-\sum_{j=h_{0}}^{h_{1}}d_{i_{j}}}
\cdots
\sum_{c_{s}=0}^{i_{h_{s}}-i_{h_{s-1}}-1-\sum_{j=h_{s-1}}^{h_{s}}d_{i_{j}}}
\no\\
&&\quad\quad\quad\quad\quad\quad\quad\quad
\prod_{n=1}^{s} a_{d_{i_{h_{n-1}+1}}\cdots d_{i_{h_{n}-1}}c_{n}}
(i_{h_{n-1}+1}-i_{h_{n-1}},\cdots,i_{h_{n}}-i_{h_{n}-1})\no\\
&&\sum_{l_{r_{1}}=1}^{c_{r_{1}+1}-c_{r_{1}}+
i_{h_{r_{1}}}-i_{h_{r_{1}-1}}-\sum_{j=h_{r_{1}-1}}^{h_{r_{1}}}d_{i_{j}}}
\cdots
\sum_{l_{r_{t}}=1}^{c_{r_{t}+1}-c_{r_{t}}+
i_{h_{r_{t}}}-i_{h_{r_{t}-1}}-\sum_{j=h_{r_{t}-1}}^{h_{r_{t}}}
d_{i_{j}}}\no\\
&&\sum_{j_{s-t}=0}^{N-(N-k+1)d}\sum_{j_{s-t-1}=0}^{j_{s-t}} 
\cdots\sum_{j_{1}=0}^{j_{2}}
\biggl(\prod_{n=1}^{s-t}(1+(d-i_{h_{p_{n-1}}})w)^{j_{n}-j_{n-1}}\biggr)\cdot
\biggl(\prod_{n=1}^{t}(1+(d-i_{h_{r_{n}}})w)^{-l_{r_{n}}}\biggr)\no\\
&&\prod_{n=1}^{s-t}\prod_{m=p_{n-1}+1}^{p_{n}}
\biggl((1+(d-i_{h_{m-1}})w)^{c_{m}-c_{m-1}+
i_{h_{m-1}}-i_{h_{m-2}}-\sum_{j=h_{m-2}}^{h_{m-1}}d_{i_{j}}}\times\no\\
&&\times
\prod_{a=h_{m-1}+1}^{h_{m}}
L^{N+1,k,i_{a}-i_{a-1}}_{j_{n}-\sum_{j=1}^{a-1}
d_{i_{j}}-\sum_{r_{j}<m}l_{r_{j}}+i_{a-1}(N-k+1)}\biggr).\no\\
\label{genebound}
\end{eqnarray}
Now, what we have to show is that the bulk part coming from the ordered 
partition $0=h_{0}<h_{1}<\cdots<h_{s-1}<h_{s}=l$ cancels with the 
boundary parts coming from the ordered partition 
$0=q_{0}<q_{1}<\cdots<q_{t-1}<q_{t}=l,\;\;(\{h_{0},h_{1},\cdots,h_{s-1},h_{s}\}
\subset \{q_{0},q_{1},\cdots,q_{t-1},q_{t}\})$. 
Before general discussion on cancellation of these terms, we carry out 
computations for some lower $l$'s as warming-up's. \\
\\
$0=i_{0}<i_{1}=d$ sector:\\
In this case, there are no boundary contributions and the bulk part is 
given by, 
\begin{eqnarray}
&&\sum_{c_{1}=0}^{d-1}\sum_{j=-c_{1}}^{N-1-(N-k)d-c_{1}}a_{c_{1}}
(d)(1+dw)^{j+c_{1}}
L_{j}^{N+1,k,d}\no\\
&&=\sum_{c_{1}=0}^{d-1}\sum_{j=0}^{N-(N-k+1)d}a_{c_{1}}(d)(1+dw)^{j+c_{1}}
L_{j}^{N+1,k,d}=\prod_{j=1}^{d-1}(1+jw)\cdot
\sum_{j=0}^{N-(N-k+1)d}(1+dw)^{j}
L_{j}^{N+1,k,d}\no\\ 
\label{relver}
\end{eqnarray}
where we used (\ref{wow}) and the fact that $L_{j}^{N+1,k,d}=0$ unless
$0\leq j\leq N-(N+1-k)d$.\\
\\
$0=i_{0}<i_{1}<i_{2}=d$ sector:\\
In this sector, the summand coming from $0=h_{0}<h_{1}=1<h_{2}=2$ is 
separated into one bulk contribution and one boundary contribution 
corresponding to $0=h_{0}<h_{2}=2$:
\begin{eqnarray}
&&-\prod_{j=1}^{d-1}(1+jw)\cdot\sum_{j_{2}=0}^{N-(N-k+1)d}
\sum_{j_{1}=0}^{j_{2}}(1+dw)^{j_{1}}(1+(d-i_{1})w)^{j_{2}-j_{1}}
L_{j_{1}+(N-k+1)i_{1}}^{N+1,k,i_{1}}L_{j_{2}+(N-k+1)i_{2}}^{N+1,k,d-i_{1}}\no\\
&&-\sum_{c_{1}=0}^{i_{1}-1}\sum_{c_{2}=0}^{d-i_{1}-1}
\sum_{j_{1}=0}^{N-(N-k+1)d}(1+dw)^{j_{1}}\sum_{l_{1}=1}^{c_{2}-c_{1}+i_{1}}
a_{c_{1}}(i_{1})a_{c_{2}}(d-i_{1})(1+(d-i_{1})w)^{-l_{1}+c_{2}-c_{1}+i_{1}}
(1+dw)^{c_{1}}\times\no\\
&&\times L_{j_{1}}^{N+1,k,i_{1}}
L_{j_{1}-l_{1}+(N-k+1)d_{1}}^{N+1,k,d-i_{1}}.
\label{reldub}
\end{eqnarray}
On the other hand, we have to prove that the second summand of 
the r.h.s of (\ref{reldub}) cancels with the summand (\ref{addend}) 
coming from $0=h_{0}<h_{1}=2$,
\begin{eqnarray}
&&\sum_{j_{1}=0}^{N-d(N-k+1)}\sum_{d_{i_{1}}=1}^{d-1}
\sum_{c_{1}=0}^{d-1-d_{i_{1}}}(1+dw)^{j_{1}}a_{d_{i_{1}}c_{1}}(i_{1},d-i_{1})
(1+dw)^{c_{1}} L_{j_{1}}^{N+1,k,i_{1}}
L_{j_{1}-d_{i_{1}}+(N-k+1)d_{1}}^{N+1,k,d-i_{1}},\no\\
\label{hoo}
\end{eqnarray} 
where 
\begin{eqnarray}
&&\sum_{j=0}^{d-1-d_{i_{1}}}a_{d_{i_{1}}j}
(i_{1},d-i_{1})x^{j}y^{d-1-d_{i_{1}}-j}\no\\
&&=\frac{1}{2\pi\sqrt{-1}}\int_{D(0,\frac{i_{1}y+(d-i_{1})x}{d})}du_{i_{1}}
\frac{1}{(u_{i_{1}}-\frac{i_{1}y+(d-i_{1})x}{d})}
\frac{f_{(i_{1})}(x,u_{i_{1}})\cdot
f_{(d-i_{1})}(u_{i_{1}},y)}{u_{i_{1}}^{d_{i_{1}}-1}}\no\\
&&=f_{(i_{1},d-i_{1})}^{d_{i_{1}}}(x,y).
\label{d1d2}
\end{eqnarray}  
Since $-l_{1}+c_{2}-c_{1}+i_{1}\geq 0$ in (\ref{reldub}),
 the assertion of the theorem in this sector 
reduces to the following polynomial 
identity in this sector:
\begin{equation}
\frac{f_{(i_{1})}(x,u)f_{(d-i_{1})}(u,y)}{u^{d_{i_{1}}-1}}|_{\deg(u)\geq 0,\;
u=\frac{(d-i_{1})x+i_{1}y}{d}}
=f_{(i_{1},d-i_{1})}^{d_{i_{1}}}(x,y), 
\label{comcom}
\end{equation}
where $\frac{f_{(i_{1})}(x,u)f_{(d-i_{1})}(u,y)}{u^{d_{i_{1}}-1}}
|_{\deg(u)\geq 0}$ means the operation of picking up monomials, whose 
degree in $u$ is non-negative, from  $\frac{f_{(i_{1})}(x,u)f_{(d-i_{1})}(u,y)}{u^{d_{i_{1}}-1}}$.  
Now, we prove the above equality. Using the residue integral in $v$-plane, 
we have,
\begin{eqnarray}
&&\frac{f_{(i_{1})}(x,u)f_{(d-i_{1})}(u,y)}{u^{d_{i_{1}}-1}}
|_{\deg(u)\geq 0}
=\frac{1}{2\pi\sqrt{-1}}\int_{D_{v}(0)}dv
\frac{f_{(i_{1})}(x,v)f_{(d-i_{1})}(v,y)}{v^{d_{i_{1}}-1}}
\frac{1}{v}\sum_{n=0}^{\infty}
(\frac{u}{v})^{n}\no\\
&=&\frac{1}{2\pi\sqrt{-1}}\int_{D_{v}(0,u)}dv
\frac{f_{(i_{1})}(x,v)f_{(d-i_{1})}(v,y)}{v^{d_{i_{1}}-1}}\frac{1}{v-u}.
\end{eqnarray} 
Therefore, we can rewrite the l. h. s. of (\ref{comcom}) as follows:
\begin{eqnarray}
&&\frac{f_{(i_{1})}(x,u)f_{(d-i_{1})}(u,y)}{u^{d_{i_{1}}-1}}|_{deg(u)\geq 0,\;
u=\frac{(d-i_{1})x+i_{1}y}{d}}\no\\
&&=(\frac{1}{2\pi\sqrt{-1}})
\int_{D_{u}(0,\frac{(d-i_{1})x+i_{1}y}{d})}du
\frac{1}{(u-\frac{(d-i_{1})x+i_{1}y}{d})}
\frac{f_{(i_{1})}(x,u)f_{(d-i_{1})}(u,y)}{u^{d_{i_{1}}-1}}.\no\\
\end{eqnarray}
But the last line is nothing but the definition of 
$f_{(i_{1},d-i_{1})}^{d_{i_{1}}}(x,y)$.\\
\\
$0=i_{0}<i_{1}<i_{2}<i_{3}=d$ sector:\\
In this case, we have four choices of partitions: 
\begin{eqnarray}
&&0=h_{0}<1=h_{1}<2=h_{2}<h_{3}=3,\no\\
&&0=h_{0}<h_{1}=1<h_{2}=3,\no\\
&&0=h_{0}<h_{1}=2<h_{2}=3,\no\\
&&0=h_{0}<h_{1}=3.
\end{eqnarray}
The summand in (\ref{addend}) coming from the ordered partition
$0=h_{0}<1=h_{1}<2=h_{2}<h_{3}=3$ is decomposed as follows:
\begin{eqnarray}
&&\prod_{j=1}^{d-1}(1+jw)\cdot
\sum_{j_{3}=0}^{N-(N-k+1)d}\sum_{j_{2}=0}^{j_{3}}\sum_{j_{1}=0}^{j_{2}}
(1+dw)^{j_{1}}(1+(d-i_{1})w)^{j_{2}-j_{1}}
(1+(d-i_{2})w)^{j_{3}-j_{2}}\times\no\\
&&\times L_{j_{1}}^{N+1,k,i_{1}} L_{j_{2}+i_{1}(N-k+1)}^{N+1,k,i_{2}-i_{1}}
 L_{j_{3}+i_{2}(N-k+1)}^{N+1,k,d-i_{2}}\no\\
&&+\sum_{c_{1}=0}^{i_{1}-1}\sum_{c_{2}=0}^{i_{2}-i_{1}-1}
\sum_{c_{3}=0}^{i_{3}-i_{2}-1}a_{c_{1}}(i_{1})a_{c_{2}}(i_{2}-i_{1})
a_{c_{3}}(d-i_{2})\times\no\\
&&\times\sum_{l_{1}=1}^{c_{2}-c_{1}+i_{1}}
\sum_{j_{3}=0}^{N-(N-k+1)d}\sum_{j_{1}=0}^{j_{3}}
(1+dw)^{j_{1}+c_{1}}(1+(d-i_{1})w)^{-l_{1}+c_{2}-c_{1}+i_{1}}\times\no\\
&&\times(1+(d-i_{2})w)^{j_{3}-j_{1}+c_{3}-c_{2}+i_{2}-i_{1}}
 L_{j_{1}}^{N+1,k,i_{1}} 
L_{j_{1}-l_{1}+i_{1}(N-k+1)}^{N+1,k,i_{2}-i_{1}}
 L_{j_{3}-l_{1}+i_{2}(N-k+1)}^{N+1,k,d-i_{2}}\no\\
&&+\sum_{c_{1}=0}^{i_{1}-1}\sum_{c_{2}=0}^{i_{2}-i_{1}-1}
\sum_{c_{3}=0}^{i_{3}-i_{2}-1}a_{c_{1}}(i_{1})a_{c_{2}}(i_{2}-i_{1})
a_{c_{3}}(d-i_{2})\times\no\\
&&\times\sum_{l_{2}=1}^{c_{3}-c_{2}+i_{2}-i_{1}}
\sum_{j_{2}=0}^{N-(N-k+1)d}\sum_{j_{1}=0}^{j_{2}}
(1+dw)^{j_{1}+c_{1}}(1+(d-i_{1})w)^{j_{2}-j_{1}+c_{2}-c_{1}+i_{1}}\times\no\\
&&\times(1+(d-i_{2})w)^{-l_{2}+c_{3}-c_{2}+i_{2}-i_{1}}
L_{j_{1}}^{N+1,k,i_{1}} L_{j_{2}+i_{1}(N-k+1)}^{N+1,k,i_{2}-i_{1}}
 L_{j_{2}-l_{2}+i_{2}(N-k+1)}^{N+1,k,d-i_{2}}\no\\
&&+\sum_{c_{1}=0}^{i_{1}-1}\sum_{c_{2}=0}^{i_{2}-i_{1}-1}
\sum_{c_{3}=0}^{i_{3}-i_{2}-1}a_{c_{1}}(i_{1})a_{c_{2}}(i_{2}-i_{1})
a_{c_{3}}(d-i_{2})\times\no\\
&&\times\sum_{l_{1}=1}^{c_{2}-c_{1}+i_{1}}
\sum_{l_{2}=1}^{c_{3}-c_{2}+i_{2}-i_{1}}\sum_{j_{1}=0}^{N-(N-k+1)d}
(1+dw)^{j_{1}+c_{1}}(1+(d-i_{1})w)^{-l_{1}+c_{2}-c_{1}+i_{1}}\times\no\\
&&\times(1+(d-i_{2})w)^{-l_{2}+c_{3}-c_{2}+i_{2}-i_{1}}
 L_{j_{1}}^{N+1,k,i_{1}} 
L_{j_{1}-l_{1}+i_{1}(N-k+1)}^{N+1,k,i_{2}-i_{1}}
 L_{j_{1}-l_{1}-l_{2}+i_{2}(N-k+1)}^{N+1,k,d-i_{2}}.
\label{123}
\end{eqnarray}
Note that the first summand, the second and the third ones, and the 
last one correspond to the bulk part, the first boundary parts, and the second 
boundary part respectively. 
Next, we decompose the summand coming from $0=h_{0}<h_{1}=1<h_{2}=3$,
\begin{eqnarray}
&&-\sum_{d_{2}=1}^{i_{2}-i_{1}-1}\sum_{c_{1}=0}^{i_{1}-1}
\sum_{c_{2}=0}^{i_{3}-i_{1}-1-d_{2}}a_{c_{1}}(i_{1})a_{d_{2}c_{2}}
(i_{2}-i_{1},d-i_{2})\times\no\\
&&\times\sum_{j_{2}=0}^{N-(N-k+1)d}\sum_{j_{1}=0}^{j_{2}}
(1+dw)^{j_{1}+c_{1}}(1+(d-i_{1})w)^{j_{2}-j_{1}+c_{2}-c_{1}+i_{1}}\times\no\\
&&\times L_{j_{1}}^{N+1,k,i_{1}} L_{j_{2}+i_{1}(N-k+1)}^{N+1,k,i_{2}-i_{1}}
 L_{j_{2}-d_{2}+i_{2}(N-k+1)}^{N+1,k,d-i_{2}}\no\\
&&-\sum_{d_{2}=1}^{i_{2}-i_{1}-1}\sum_{c_{1}=0}^{i_{1}-1}
\sum_{c_{2}=0}^{i_{3}-i_{1}-1-d_{2}}a_{c_{1}}(i_{1})a_{d_{2}c_{2}}
(i_{2}-i_{1},d-i_{2})\times\no\\
&&\times\sum_{l_{1}=1}^{c_{2}-c_{1}+i_{1}}\sum_{j_{1}=0}^{N-(N-k+1)d}
(1+dw)^{j_{1}+c_{1}}(1+(d-i_{1})w)^{-l_{1}+c_{2}-c_{1}+i_{1}}\times\no\\
&&\times L_{j_{1}}^{N+1,k,i_{1}} 
L_{j_{1}-l_{1}+i_{1}(N-k+1)}^{N+1,k,i_{2}-i_{1}}
 L_{j_{1}-l_{1}-d_{2}+i_{2}(N-k+1)}^{N+1,k,d-i_{2}},\no\\
\label{13}
\end{eqnarray}
and the one from $0=h_{0}<h_{1}=2<h_{2}=3$,
\begin{eqnarray}
&&-\sum_{d_{1}=1}^{i_{2}-1}\sum_{c_{1}=0}^{i_{2}-1-d_{1}}
\sum_{c_{2}=0}^{i_{3}-i_{2}-1}a_{d_{1}c_{1}}(i_{1},i_{2}-i_{1})a_{c_{2}}
(d-i_{2})\times\no\\
&&\times\sum_{j_{2}=0}^{N-(N-k+1)d}\sum_{j_{1}=0}^{j_{2}}
(1+dw)^{j_{1}+c_{1}}(1+(d-i_{2})w)^{j_{2}-j_{1}+c_{2}-c_{1}+i_{2}-d_{1}}
\times\no\\
&&\times L_{j_{1}}^{N+1,k,i_{1}} 
L_{j_{1}-d_{1}+i_{1}(N-k+1)}^{N+1,k,i_{2}-i_{1}}
 L_{j_{2}-d_{1}+i_{2}(N-k+1)}^{N+1,k,d-i_{2}}\no\\
&&-\sum_{d_{1}=1}^{i_{2}-1}\sum_{c_{1}=0}^{i_{2}-1-d_{1}}
\sum_{c_{2}=0}^{i_{3}-i_{2}-1}a_{d_{1}c_{1}}(i_{1},i_{2}-i_{1})a_{c_{2}}
(d-i_{2})\times\no\\
&&\times\sum_{l_{1}=1}^{c_{2}-c_{1}+i_{2}-d_{1}}
\sum_{j_{1}=0}^{N-(N-k+1)d}
(1+dw)^{j_{1}+c_{1}}(1+(d-i_{2})w)^{-l_{1}+c_{2}-c_{1}+i_{2}-d_{1}}\times\no\\
&&\times L_{j_{1}}^{N+1,k,i_{1}} 
L_{j_{1}-d_{1}+i_{1}(N-k+1)}^{N+1,k,i_{2}-i_{1}}
 L_{j_{1}-d_{1}-l_{1}+i_{2}(N-k+1)}^{N+1,k,d-i_{2}}.
\label{23}
\end{eqnarray}
The summand coming from $0=h_{0}<h_{3}=3$ is given by,
\begin{eqnarray}
&&\sum_{j_{1}=0}^{N-d(N-k+1)}\sum_{d_{1}=1}^{d-1}
\sum_{d_{2}=1}^{d-1-d_{1}}
\sum_{c_{1}=0}^{d-1-d_{1}-d_{2}}a_{d_{1}d_{2}c_{1}}
(i_{1},i_{2}-i_{1},d-i_{2})(1+dw)^{j_{1}+c_{1}}
L_{j_{1}}^{N+1,k,i_{1}}\times\no\\
&&\times L_{j_{1}-d_{1}+i_{1}(N-k+1)}^{N+1,k,i_{2}-i_{1}}
L_{j_{1}-d_{1}-d_{2}+i_{2}(N-k+1)}^{N+1,k,d-i_{2}}.
\label{eq6}
\end{eqnarray}
With these results, we can easily see that the second (resp. third) 
summand in (\ref{123})
cancels with the first summand in (\ref{23}) (resp. (\ref{13})) due to 
the identity proved in the $l=2$ case. Therefore, the new identity we 
have to prove comes from the cancellation of the fourth summand in (\ref{123})
, the second summand in (\ref{13}) and in (\ref{23}), and (\ref{eq6}).
This can be translated into the polynomial equality:
\begin{eqnarray}
&&\sum_{c_{1}=0}^{i_{1}-1}\sum_{c_{2}=0}^{i_{2}-i_{1}-1}
\sum_{c_{3}=0}^{i_{3}-i_{2}-1}a_{c_{1}}(i_{1})a_{c_{2}}(i_{2}-i_{1})
a_{c_{3}}(d-i_{2})\times\no\\
&&\times\sum_{l_{1}=1}^{c_{2}-c_{1}+i_{1}}
\sum_{l_{2}=1}^{c_{3}-c_{2}+i_{2}-i_{1}}
(1+dw)^{c_{1}}(1+(d-i_{1})w)^{-l_{1}+c_{2}-c_{1}+i_{1}}
(1+(d-i_{2})w)^{-l_{2}+c_{3}-c_{2}+i_{2}-i_{1}}\no\\
&&-\sum_{d_{2}=0}^{i_{2}-i_{1}-1}\sum_{c_{1}=0}^{i_{1}-1}
\sum_{c_{2}=0}^{i_{3}-i_{1}-1-d_{2}}a_{c_{1}}(i_{1})a_{d_{2}c_{2}}
(i_{2}-i_{1},d-i_{2})\sum_{l_{1}=1}^{c_{2}-c_{1}+i_{1}}
(1+dw)^{c_{1}}(1+(d-i_{1})w)^{-l_{1}+c_{2}-c_{1}+i_{1}}\no\\
&&-\sum_{d_{1}=0}^{i_{2}-1}\sum_{c_{1}=0}^{i_{1}-1-d_{1}}
\sum_{c_{2}=0}^{i_{3}-i_{2}-1}a_{d_{1}c_{1}}(i_{1},i_{2}-i_{1})a_{c_{2}}
(d-i_{2})\sum_{l_{1}=1}^{c_{2}-c_{1}+i_{2}-d_{1}}
(1+dw)^{c_{1}}(1+(d-i_{2})w)^{-l_{1}+c_{2}-c_{1}+i_{2}-d_{1}}\no\\
&&+\sum_{d_{1}=1}^{d-1}
\sum_{d_{2}=1}^{d-1-d_{1}}
\sum_{c_{1}=0}^{d-1-d_{1}-d_{2}}a_{d_{1}d_{2}c_{1}}
(i_{1},i_{2}-i_{1},d-i_{2})(1+dw)^{c_{1}}=0.
\label{ni}
\end{eqnarray}
We can rewrite the above condition in a more compact form,
\begin{eqnarray}
&&\biggl(\frac{f_{(i_{1})}(x,u)f^{n}_{(i_{2}-i_{1},d-i_{2})}(u,y)}{u^{m-1}}
|_{\deg(u)\geq 0}+\frac{f^{m}_{(i_{1},i_{2}-i_{1})}(x,v)
f_{(d-i_{2})}(v,y)}{v^{n-1}}
|_{\deg(v)\geq 0}\no\\
&&-\frac{f_{(i_{1})}(x,u)f_{(i_{2}-i_{1})}(u,v)f_{(d-i_{2})}(v,y)}
{u^{m-1}v^{n-1}}
|_{\deg(u)\geq 0,\;\deg(v)\geq 0}
\biggr)
|_{u=\frac{i_{1}y+(d-i_{1})x}{d},
v=\frac{(i_{2})y+(d-i_{2})x}{d}}\no\\ 
&&=f^{mn}_{(i_{1},i_{2}-i_{1},d-i_{2})}(x,y).
\label{3}
\end{eqnarray}
On the other hand, the definition of 
$f^{mn}_{(i_{1},i_{2}-i_{1},d-i_{2})}(x,y)$
tells us, 
\begin{eqnarray}
&&f_{(i_{1},i_{2}-i_{1},d-i_{2})}^{mn}(x,y):=\no\\
&&(\frac{1}{2\pi\sqrt{-1}})^{2}\int_{D_{u}}\int_{D_{v}}
\frac{d\cdot (i_{2}-i_{1})\cdot(du dv)}{(i_{2}u-i_{1}v-(i_{2}-i_{1})x)
((d-i_{1})v-(d-i_{2})u-(i_{2}-i_{1})y)}\times\no\\
&&\times\frac{f_{(i_{1})}(x,u)f_{(i_{2}-i_{1})}(u,v)
f_{(d-i_{2})}(v,y)}{u^{m-1}v^{n-1}}.\no\\
\label{sim2001}
\end{eqnarray} 
Hence what remains to show in the $l=3$ case is the following equality:
\begin{eqnarray}
&&\biggl(\frac{f_{(i_{1})}(x,u)f^{n}_{(i_{2}-i_{1},d-i_{2})}(u,y)}{u^{m-1}}
|_{\deg(u)\geq 0}+\frac{f^{m}_{(i_{1},i_{2}-i_{1})}(x,v)
f_{(d-i_{2})}(v,y)}{v^{n-1}}
|_{\deg(v)\geq 0}\no\\
&&-\frac{f_{(i_{1})}(x,u)f_{(i_{2}-i_{1})}(u,v)f_{(d-i_{2})}(v,y)}
{u^{m-1}v^{n-1}}
|_{\deg(u)\geq 0,\;\deg(v)\geq 0}
\biggr)
|_{u=\frac{i_{1}y+(d-i_{1})x}{d},
v=\frac{(i_{2})y+(d-i_{2})x}{d}}\no\\ 
&&=(\frac{1}{2\pi\sqrt{-1}})^{2}\int_{D_{u}}\int_{D_{v}}
\frac{d\cdot (i_{2}-i_{1})\cdot(du dv)}{(i_{2}u-i_{1}v-(i_{2}-i_{1})x)
((d-i_{1})v-(d-i_{2})u-(i_{2}-i_{1})y)}\times\no\\
&&\times\frac{f_{(i_{1})}(x,u)f_{(i_{2}-i_{1})}(u,v)
f_{(d-i_{2})}(v,y)}{u^{m-1}v^{n-1}}.
\label{3e}
\end{eqnarray}
First, we consider the following part:
\begin{eqnarray}
&&(\frac{1}{2\pi\sqrt{-1}})^{2}\int_{D_{u}}\int_{D_{v}}
\biggl(\frac{f_{(i_{1})}(x,u)f_{(i_{2}-i_{1})}(u,v)
f_{(d-i_{2})}(v,y)}{u^{m-1}v^{n-1}}
|_{\deg(u)\leq -1, \deg(v)\leq -1}\biggr)\times\no\\
&&\times\frac{d\cdot (i_{2}-i_{1})\cdot(du dv)}{((i_{2})u-i_{1}v-
(i_{2}-i_{1})x)
((d-i_{1})v-(d-i_{2})u-(i_{2}-i_{1})y)}.
\end{eqnarray} 
But we can easily see with some computation,    
\begin{eqnarray}
&&(\frac{1}{2\pi\sqrt{-1}})^{2}\int_{D_{u}}\int_{D_{v}}
\frac{1}{u^{k}v^{l}}\cdot
\frac{d\cdot d_{2}\cdot(du dv)}{((i_{2})u-i_{1}v-(i_{2}-i_{1})x)
((d-i_{1})v-(d-i_{2})u-(i_{2}-i_{1})y)}\no\\
&&=\frac{1}{2\pi\sqrt{-1}}\int_{D_{u}}
\frac{du}{u^{k}}\cdot\frac{(i_{1})^{l}}{((i_{2})u-(i_{2}-i_{1})x)^{l}}
\cdot \frac{1}{(u-\frac{i_{1}y+(d-i_{1})x}{d})}\no\\
&&=0,
\end{eqnarray}
where $k,l \geq 1$. The last equality follows from the fact that $D_{u}$ 
goes  around all the poles of the integrand. Hence we have 
\begin{eqnarray}
&&(\frac{1}{2\pi\sqrt{-1}})^{2}\int_{D_{u}}\int_{D_{v}}
\biggl(\frac{f_{(i_{1})}(x,u)f_{(i_{2}-i_{1})}(u,v)
f_{(d-i_{2})}(v,y)}{u^{m-1}v^{n-1}}
|_{\deg(u)\leq -1, \deg(v)\leq -1}\biggr)\times\no\\
&&\times\frac{d\cdot (i_{2}-i_{1})\cdot(du dv)}{((i_{2})u-i_{1}v-
(i_{2}-i_{1})x)
((d-i_{1})v-(d-i_{2})u-(i_{2}-i_{1})y)}\no\\
&&=0.
\label{van}
\end{eqnarray}
Using (\ref{van}), we can rewrite the r.h.s. of (\ref{3e}) as follows,
\begin{eqnarray}
&&(\frac{1}{2\pi\sqrt{-1}})^{2}\int_{D_{u}}\int_{D_{v}}
\frac{d\cdot (i_{2}-i_{1})\cdot(du dv)}{(i_{2}u-i_{1}v-(i_{2}-i_{1})x)
((d-i_{1})v-(d-i_{2})u-(i_{2}-i_{1})y)}\times\no\\
&&\times\frac{f_{(i_{1})}(x,u)f_{(i_{2}-i_{1})}(u,v)
f_{(d-i_{2})}(v,y)}{u^{m-1}v^{n-1}}\no\\
&&=(\frac{1}{2\pi\sqrt{-1}})^{2}\int_{D_{u}}\int_{D_{v}}
\frac{d\cdot (i_{2}-i_{1})\cdot(du dv)}{(i_{2}u-i_{1}v-(i_{2}-i_{1})x)
((d-i_{1})v-(d-i_{2})u-(i_{2}-i_{1})y)}\times\no\\
&&\times\frac{f_{(i_{1})}(x,u)f_{(i_{2}-i_{1})}(u,v)
f_{(d-i_{2})}(v,y)}{u^{m-1}v^{n-1}}\no\\
&&-(\frac{1}{2\pi\sqrt{-1}})^{2}\int_{D_{u}}\int_{D_{v}}
\frac{d\cdot (i_{2}-i_{1})\cdot(du dv)}{(i_{2}u-i_{1}v-(i_{2}-i_{1})x)
((d-i_{1})v-(d-i_{2})u-(i_{2}-i_{1})y)}\times\no\\
&&\times\biggl(\frac{f_{(i_{1})}(x,u)f_{(i_{2}-i_{1})}(u,v)
f_{(d-i_{2})}(v,y)}{u^{m-1}v^{n-1}}|_{\deg(u)\leq -1,\;\;g(v)\leq -1}
\biggr)\no\\
&&=(\frac{1}{2\pi\sqrt{-1}})^{2}\int_{D_{u}}\int_{D_{v}}
\frac{d\cdot (i_{2}-i_{1})\cdot(du dv)}{(i_{2}u-i_{1}v-(i_{2}-i_{1})x)
((d-i_{1})v-(d-i_{2})u-(i_{2}-i_{1})y)}\times\no\\
&&\times\biggl(\frac{f_{(i_{1})}(x,u)f_{(i_{2}-i_{1})}(u,v)
f_{(d-i_{2})}(v,y)}{u^{m-1}v^{n-1}}|_{\deg(u)\geq 0}
\biggr)\no\\
&&+(\frac{1}{2\pi\sqrt{-1}})^{2}\int_{D_{u}}\int_{D_{v}}
\frac{d\cdot (i_{2}-i_{1})\cdot(du dv)}{(i_{2}u-i_{1}v-(i_{2}-i_{1})x)
((d-i_{1})v-(d-i_{2})u-(i_{2}-i_{1})y)}\times\no\\
&&\times\biggl(\frac{f_{(i_{1})}(x,u)f_{(i_{2}-i_{1})}(u,v)
f_{(d-i_{2})}(v,y)}{u^{m-1}v^{n-1}}|_{\deg(v)\geq 0}
\biggr)\no\\
&&-(\frac{1}{2\pi\sqrt{-1}})^{2}\int_{D_{u}}\int_{D_{v}}
\frac{d\cdot (i_{2}-i_{1})\cdot(du dv)}{(i_{2}u-i_{1}v-(i_{2}-i_{1})x)
((d-i_{1})v-(d-i_{2})u-(i_{2}-i_{1})y)}\times\no\\
&&\times\biggl(\frac{f_{(i_{1})}(x,u)f_{(i_{2}-i_{1})}(u,v)
f_{(d-i_{2})}(v,y)}{u^{m-1}v^{n-1}}|_{\deg(v)\geq 0,\;\;\deg(u)\geq 0}
\biggr).
\label{inc}
\end{eqnarray}
At this stage, we look at the first integral of the last line of (\ref{inc}).
Due to the condition $\deg(u)\geq 0$, $u$ variable has only one pole at 
$u=\frac{i_{1}v+(i_{2}-i_{1})x}{i_{1}}$. And if we integrate out the 
$u$ variable first, the integral turns into,
\begin{equation}
(\frac{1}{2\pi\sqrt{-1}})\int_{D_{v}}
\frac{dv}{(v-\frac{i_{2}y+(d-i_{2})x}{d})}
\biggl(\frac{f^{m}_{(i_{1},i_{2}-i_{1})}(x,v)f_{(d-i_{2})}(v,y)}
{v^{n-1}}\biggr).
\end{equation}
This is nothing but the second term in the l.h.s. of (\ref{3e}). Using the 
same operation, we can show that the second and the third integrals in the 
last line of (\ref{inc}) equal the first and the third terms in the l.h.s.
of (\ref{3e}). Thus, the proof of $l=3$ case is completed.   

With these preparation, we turn into the general proof of the theorem. 
In this case, we have to consider the integral,
\begin{eqnarray}
&&f_{(i_{1}-i_{0},\cdots,i_{l}-i_{l-1})}^{d_{i_{1}}d_{i_{2}}\cdots d_{i_{l-1}}}
(x,y)=\no\\
&&\frac{d}{(2\pi\sqrt{-1})^{l-1}}\int_{D_{i_{1}}}\cdots
\int_{D_{i_{l-1}}}\biggl(\prod_{j=1}^{l}
f_{(i_{j}-i_{j-1})}(u_{i_{j-1}},u_{i_{j}})
\prod_{j=1}^{l-1}\frac{1}{u_{i_{j}}^{d_{i_{j}}-1}}\biggr)\times\no\\
&&\times\prod_{j=1}^{l-1}\frac{du_{i_{j}}}{((i_{j+1}-i_{j-1})u_{i_{j}}-
(i_{j}-i_{j-1})u_{i_{j+1}}-(i_{j+1}-i_{j})u_{i_{j-1}})}
\prod_{j=1}^{l-2}(i_{j+1}-i_{j}).
\end{eqnarray}
For convenience of space, we introduce the definition:
\begin{defi}
Let $\alpha_{j}(x,y)\;\;(j=1,2,\cdots,l)$ be a 
homogeneous polynomial in $x$ and $y$. We define 
two types of $l$-product, which are both non-commutative 
and non-associative as follows,
\begin{eqnarray}
&&(\alpha_{1}\circ\alpha_{2}\circ\cdots\circ\alpha_{l})(u_{i_{0}},u_{i_{l}})
:=\no\\
&&\frac{(i_{l}-i_{0})}{(2\pi\sqrt{-1})^{l-1}}\int_{D_{i_{1}}}\cdots
\int_{D_{i_{l-1}}}\biggl(\prod_{j=1}^{l}
\alpha_{j}(u_{i_{j-1}},u_{i_{j}})
\prod_{j=1}^{l-1}\frac{1}{u_{i_{j}}^{d_{i_{j}}-1}}\biggr)\times\no\\
&&\times\prod_{j=1}^{l-1}\frac{du_{i_{j}}}{((i_{j+1}-i_{j-1})u_{i_{j}}-
(i_{j}-i_{j-1})u_{i_{j+1}}-(i_{j+1}-i_{j})u_{i_{j-1}})}
\prod_{j=1}^{l-2}(i_{j+1}-i_{j}),\\
&&(\alpha_{1}*\alpha_{2}*\cdots*\alpha_{l})(u_{i_{0}},u_{i_{l}}):=\no\\
&&\biggl(\prod_{j=1}^{l}
\alpha_{j}(u_{i_{j-1}},u_{i_{j}})
\prod_{j=1}^{l-1}\frac{1}{u_{i_{j}}^{d_{i_{j}}-1}}\biggr)
|_{\deg(u_{i_{1}})\geq 0,\;\cdots,\;\deg(u_{i_{l-1}})\geq 0,
\;\;u_{i_{j}}=\frac{(i_{j}-i_{0})u_{i_{l}}+(i_{l}-i_{j})u_{i_{0}}}
{(i_{l}-i_{0})}}.\no\\
\end{eqnarray}
\end{defi}
In the same way as the $l=3$ cases, we can show the following two lemmas.
\begin{lem}
\begin{eqnarray}
&&\frac{(i_{l}-i_{0})}{(2\pi\sqrt{-1})^{l-1}}\int_{D_{i_{1}}}\cdots
\int_{D_{i_{l-1}}}\biggl(\prod_{j=1}^{l}
\alpha_{j}(u_{i_{j-1}},u_{i_{j}})
\prod_{j=1}^{l-1}\frac{1}{u_{i_{j}}^{d_{i_{j}}-1}}|_{
\deg(u_{i_{1}})\leq -1,\cdots,\deg(u_{i_{l-1}})\leq -1}\biggr)\times\no\\
&&\times\prod_{j=1}^{l-1}\frac{du_{i_{j}}}{((i_{j+1}-i_{j-1})u_{i_{j}}-
(i_{j}-i_{j-1})u_{i_{j+1}}-(i_{j+1}-i_{j})u_{i_{j-1}})}
\prod_{j=1}^{l-2}(i_{j+1}-i_{j})=0.
\label{min}
\end{eqnarray}
\end{lem}
{\it proof)}
By expanding $\frac{1}{((i_{j+1}-i_{j-1})u_{i_{j}}-
(i_{j}-i_{j-1})u_{i_{j+1}}-(i_{j+1}-i_{j})u_{i_{j-1}})}$ into
\begin{eqnarray}
&&\frac{1}{(i_{j+1}-i_{j-1})u_{i_{j}}}\cdot
\biggl(1-\frac{(i_{j}-i_{j-1})u_{i_{j+1}}+(i_{j+1}-i_{j})u_{i_{j-1}}}
{(i_{j+1}-i_{j-1})u_{i_{j}}}\biggr)^{-1}\no\\
&&=\frac{1}{(i_{j+1}-i_{j-1})u_{i_{j}}}\cdot
\sum_{n=0}^{\infty}\biggl(\frac{(i_{j}-i_{j-1})u_{i_{j+1}}+
(i_{j+1}-i_{j})u_{i_{j-1}}}{(i_{j+1}-i_{j-1})u_{i_{j}}}\biggr)^{n},
\end{eqnarray}
we can reduce the integral in the l.h.s. of (\ref{min}) to (infinite) linear 
combination of the integral:
\begin{equation}
\frac{1}{(2\pi\sqrt{-1})^{l-1}}\int_{C_{i_{1}}}du_{i_{1}}\cdots
\int_{C_{i_{l-1}}}du_{i_{l-1}}\prod_{j=0}^{l}(u_{i_{j}})^{m_{i_{j}}},
\label{mono}
\end{equation} 
where the path ${C_{i_{j}}}$ goes around $u_{i_{j}}=0$. But we can easily 
see $m_{i_{0}}, m_{i_{l}}\geq 0$ and $\sum_{j=0}^{l}m_{i_{j}}\leq -l$, due 
to the condition $\deg(u_{i_{1}})\leq -1,\cdots,\deg(u_{i_{l-1}})\leq -1$.
Therefore, $\sum_{j=1}^{l-1}m_{i_{j}}\leq -l$. It follows that there exists 
$j\in\{1,2,\cdots,l-1\}$ such that $\deg(u_{i_{j}})$ is less than $-1$. 
Hence the integral (\ref{mono}) vanishes and the lemma is proved. Q.E.D. 
\begin{lem}
\begin{eqnarray}
&&\alpha_{1}\circ\cdots\circ\alpha_{l}\no\\
&&=\sum_{s=1}^{l-1}(-1)^{l-1-s}
\sum_{0=h_{0}<\cdots<h_{s}=l}(\alpha_{1}*\cdots*\alpha_{h_{1}})\circ
(\alpha_{h_{1}+1}*\cdots*\alpha_{h_{2}})\circ\cdots\circ(\alpha_{h_{s-1}+1}*
\cdots*\alpha_{h_{s}}).\no\\
\label{pro2}
\end{eqnarray}
\end{lem}
{\it proof)}
We denote by $A_{j}(F)$ the operation picking up the monomials with 
$\deg(u_{i_{j}})\geq 0$ from $F=\prod_{j=1}^{l}
\alpha_{j}(u_{i_{j-1}},u_{i_{j}})
\prod_{j=1}^{l-1}\frac{1}{u_{i_{j}}^{d_{i_{j}}-1}}$.
 Using Lemma 1 and the inclusion-exclusion 
principle, we obtain,
\begin{eqnarray}
&&(\alpha_{1}\circ\cdots\circ\alpha_{l})(u_{i_{0}},u_{i_{l}})\no\\
&&=\frac{(i_{l}-i_{0})}{(2\pi\sqrt{-1})^{l-1}}\int_{D_{i_{1}}}\cdots
\int_{D_{i_{l-1}}}\biggl(\prod_{j=1}^{l}
\alpha_{j}(u_{i_{j-1}},u_{i_{j}})
\prod_{j=1}^{l-1}\frac{1}{u_{i_{j}}^{d_{i_{j}}-1}}\biggr)\times\no\\
&&\times\prod_{j=1}^{l-1}\frac{du_{i_{j}}}{((i_{j+1}-i_{j-1})u_{i_{j}}-
(i_{j}-i_{j-1})u_{i_{j-1}}-(i_{j+1}-i_{j})u_{i_{j+1}})}
\prod_{j=1}^{l-2}(i_{j+1}-i_{j})\no\\
&&=\frac{(i_{l}-i_{0})}{(2\pi\sqrt{-1})^{l-1}}\int_{D_{i_{1}}}\cdots
\int_{D_{i_{l-1}}}\biggl((\cup_{j=1}^{l-1}A_{j})(F)
\biggr)\times\no\\
&&\times\prod_{j=1}^{l-1}\frac{du_{i_{j}}}{((i_{j+1}-i_{j-1})u_{i_{j}}-
(i_{j}-i_{j-1})u_{i_{j-1}}-(i_{j+1}-i_{j})u_{i_{j+1}})}
\prod_{j=1}^{l-2}(i_{j+1}-i_{j})\no\\
&&=\sum_{s=1}^{l-1}(-1)^{s}\sum_{1\leq h_{1}<\cdots<h_{s}\leq l-1}
\frac{(i_{l}-i_{0})}{(2\pi\sqrt{-1})^{l-1}}\int_{D_{i_{1}}}\cdots
\int_{D_{i_{l-1}}}\biggl((\cap_{j=1}^{s}A_{h_{j}})(F)
\biggr)\times\no\\
&&\times\prod_{j=1}^{l-1}\frac{du_{i_{j}}}{((i_{j+1}-i_{j-1})u_{i_{j}}-
(i_{j}-i_{j-1})u_{i_{j-1}}-(i_{j+1}-i_{j})u_{i_{j+1}})}
\prod_{j=1}^{l-2}(i_{j+1}-i_{j}).
\label{fin}
\end{eqnarray}
Integrating out $u_{i_{h_{j}}}$'s in the last line of (\ref{fin}) leads 
us to the assertion of the lemma.      Q.E.D.\\
\\ 
At this stage, we look back at the computation for the $l\leq 3$ cases and 
the structure of the boundary parts given in (\ref{genebound}). 
Then we can see that the bulk part (\ref{bulk}) 
coming from the ordered partition $0=h_{0}<h_{1}<h_{2}<
\cdots <h_{s-1}<h_{s}=l$ cancels with 
boundary parts 
obtained from the ordered partition $0=q_{0}<q_{1}<q_{2}<
\cdots <q_{t-1}<q_{t}=l,\;\;(\{h_{0},h_{1},\cdots,h_{s}\}\subset              \{q_{0},q_{1},\cdots,q_{t}\})$ if the following equality holds true.
\begin{prop}
\begin{eqnarray}
&&\sum_{s=1}^{l}(-1)^{s-1}\sum_{0=h_{0}<\cdots<h_{s}=l}
\biggl(\prod_{k=1}^{s}
f_{(i_{h_{k-1}+1}-i_{h_{k-1}},
\cdots,i_{h_{k}}-i_{h_{k}-1})}^{d_{i_{h_{k-1}+1}}\cdots 
d_{i_{h_{k}-1}}}
(u_{i_{h_{k-1}}},u_{i_{h_{k}}})\times\no\\
&&\times\prod_{j=1}^{s-1}
\frac{1}{(u_{i_{h_{j}}})^{d_{i_{h_{j}}}-1}}\biggr)
|_{\deg(u_{i_{h_{1}}})\geq 0,\deg(u_{i_{h_{2}}})\geq 0,\cdots,
\deg(u_{i_{h_{s}}})\geq 0,\;\;
u_{i_{h_{j}}}=\frac{i_{h_{j}}y+(d-i_{h_{j}})x}{d}}\no\\
&&=0.
\end{eqnarray}
\end{prop}
{\it proof)}
For the proof of the assertion of the proposition, 
it is sufficient to show the following relation,
\begin{eqnarray}
&&\alpha_{1}\circ\cdots\circ\alpha_{l}\no\\
&&=\sum_{s=2}^{l}(-1)^{s}
\sum_{0=h_{0}<\cdots<h_{s}=l}(\alpha_{1}\circ\cdots\circ\alpha_{h_{1}})*
(\alpha_{h_{1}+1}\circ\cdots\circ\alpha_{h_{2}})*\cdots*(\alpha_{h_{s-1}+1}
\circ\cdots\circ\alpha_{h_{s}}).\no\\
\label{pro1}
\end{eqnarray} 
We have to notice here that we can represent 
$\alpha_{1}\circ\cdots\circ\alpha_{l}$ in terms $*$-product by  
using iteratively (\ref{pro1}) only, or by iterative use of (\ref{pro2}).  
Therefore, to show the equivalence between (\ref{pro1}) and (\ref{pro2}),
it is enough for us to prove that both $*$-product representations
of $\alpha_{1}\circ\cdots\circ\alpha_{l}$ obtained from 
(\ref{pro1}) and (\ref{pro2}) coincide for all $l$. Since the $*$-product
$\alpha_{1}*\cdots*\alpha_{l}$ has different meaning for each $l$, we have to 
take care of the way of insertion of parenthesis $(\;\;)$ into 
$\alpha_{1}*\cdots*\alpha_{l}$. For example, we have to distinguish
$((\alpha_{1}*\alpha_{2})*\alpha_{3})*\alpha_{4}$ from
$(\alpha_{1}*\alpha_{2}*\alpha_{3})*\alpha_{4}$. Using this fact, we give 
here some symbolic discussion. First, we denote by $Q_{l}$ the set of all the
non-trivial ways of inserting parentheses into 
$\alpha_{1}*\cdots*\alpha_{l}$. Next, for $\pi_{l}\in Q_{l}$, we use the 
notation $\pi_{l}(\alpha_{1}*\cdots*\alpha_{l})$ for the result of insertion 
of parentheses. For example, 
\begin{equation}
\pi_{4}(\alpha_{1}*\alpha_{2}*\alpha_{3}*\alpha_{4})=
(\alpha_{1}*\alpha_{2})*(\alpha_{3}*\alpha_{4})
\end{equation} 
We also denote by $|\pi_{l}|$ the number of parentheses inserted by 
$\pi_{l}$.
  
With these preparation, we can easily obtain from (\ref{pro1}) the formula:
\begin{equation}
\alpha_{1}\circ\cdots\circ\alpha_{l}=\sum_{\pi_{l}\in Q_{l}}
(-1)^{(l-|\pi_{l}|)}\pi_{l}(\alpha_{1}*\cdots*\alpha_{l}),
\label{ql}
\end{equation}  
by induction of $l$. Therefore, what remains to show is that we can derive 
the formula (\ref{ql}) only by using (\ref{pro2}). We show this by induction 
of $l$. In the $l=2$ case, (\ref{pro2}) reduces to $\alpha_{1}\circ\alpha_{2}
=\alpha_{1}*\alpha_{2}$, and (\ref{ql}) trivially holds. Then we assume that 
(\ref{ql}) holds for $l=1,2,\cdots,m-1$ cases. By the assumption of induction,
it is clear that all the $\pi_{m}(\alpha_{1}*\cdots*\alpha_{m})\;\;
(\pi_{m}\in Q_{m})$ appear in the process of rewriting $\alpha_{1}
\circ\cdots\circ\alpha_{m}$ using (\ref{pro2}). Therefore, we only have to 
show that the coefficient of $\pi_{m}(\alpha_{1}*\cdots\alpha_{m})$ becomes
$(-1)^{(m-|\pi_{m}|)}$ after adding up all the contributions.    

Now, we fix one $\pi_{m}(\alpha_{1}*\cdots*\alpha_{m})$. By assumption, 
the terms coming from one term in (\ref{pro2}):
\begin{equation}  
(-1)^{l-1-s}(\alpha_{1}*\cdots*\alpha_{h_{1}})\circ
(\alpha_{h_{1}+1}*\cdots*\alpha_{h_{2}})\circ\cdots\circ(\alpha_{h_{s-1}+1}*
\cdots*\alpha_{h_{s}}),
\label{suma}
\end{equation}
are all different from each other, and we first determine the term 
(\ref{suma}) that produces $\pi_{m}(\alpha_{1}*\cdots*\alpha_{m})$.
Here, we have to notice that the terms coming from (\ref{suma}) have no
insertion of parentheses inside $(\alpha_{h_{n-1}+1}*\cdots*\alpha_{h_{n}})$.
With this observation, we remove the  
parentheses in $\pi_{m}(\alpha_{1}*\cdots*\alpha_{m})$
if they have other parentheses inside them. We denote 
by $\tilde{\pi}_{m}(\alpha_{1}*\cdots*\alpha_{m})$ the resulting term.
$\tilde{\pi}_{m}(\alpha_{1}*\cdots*\alpha_{m})$ has the following structure:
\begin{eqnarray}
&&\tilde{\pi}_{m}(\alpha_{1}*\cdots*\alpha_{m})=
\alpha_{1}*\cdots*\alpha_{k_{1}}
*(\alpha_{k_{1}+1}*\cdots*\alpha_{j_{1}})*\alpha_{j_{1}+1}*
\cdots*\alpha_{k_{2}}*(\alpha_{k_{2}+1}*\cdots*\alpha_{j_{2}})*\cdots\no\\
&&\cdots*\alpha_{k_{n}}*(\alpha_{k_{n}+1}*\cdots*\alpha_{j_{n}})
*\alpha_{j_{n}+1}*
\cdots*\alpha_{m}.
\label{piti}
\end{eqnarray} 
We determine here the terms (\ref{suma}) that produce (\ref{piti}). Since 
we cannot admit the part $(\alpha_{a_{1}}*\cdots*\alpha_{a_{1}+b})$ in
(\ref{suma}) that do not appear in (\ref{piti}), the allowed terms are
\begin{eqnarray}
&&\alpha_{1}\circ\cdots\circ\alpha_{k_{1}}
\circ(\alpha_{k_{1}+1}*\cdots*\alpha_{j_{1}})\circ\alpha_{j_{1}+1}\circ
\cdots\circ\alpha_{k_{2}}\circ(\alpha_{k_{2}+1}*\cdots*\alpha_{j_{2}})
\circ\cdots\no\\
&&\cdots\circ\alpha_{k_{n}}\circ(\alpha_{k_{n}+1}*\cdots*\alpha_{j_{n}})
\circ\alpha_{j_{n}+1}\circ
\cdots\circ\alpha_{m}.
\label{citi}
\end{eqnarray}
and the terms obtained from changing 
$\circ(\alpha_{k_{a}+1}*\cdots*\alpha_{j_{a}})\circ $ in (\ref{citi}) into 
$\circ\alpha_{k_{a}+1}\circ\cdots\circ\alpha_{j_{a}}\circ$. Here, we omit the 
sign of (\ref{citi}) for brevity. If we change all the 
$\circ(\alpha_{k_{a}+1}*\cdots*\alpha_{j_{a}})\circ $'s into 
$\circ\alpha_{k_{a}+1}\circ\cdots\circ\alpha_{j_{a}}\circ$'s, we obtain 
$\alpha_{1}\circ\cdots\circ\alpha_{l}$. Therefore, total number of the terms
(\ref{suma}) that produce $\tilde{\pi}_{m}(\alpha_{1}*\cdots*\alpha_{m})$
is $2^{n}-1=2^{|\tilde{\pi}_{m}|}-1$.
With some computation, we can see that the sign of (\ref{piti}) 
, coming from the term obtained 
from changing $h$ of the  $\circ(\alpha_{k_{a}+1}*\cdots*\alpha_{j_{a}})
\circ $'s into $\circ\alpha_{k_{a}+1}\circ\cdots\circ\alpha_{j_{a}}\circ$'s,
equals $(-1)^{m-1-h}$. And the number of such terms are given by
$n\choose h$. Therefore, the coefficient of 
$\tilde{\pi}_{m}(\alpha_{1}*\cdots*\alpha_{m})$ turns out to be,
\begin{equation}
\sum_{h=0}^{n-1}{n\choose h}(-1)^{m-1-h}=(-1)^{m-1}((1-1)^{n}-(-1)^{n})
=(-1)^{m-n}=(-1)^{(m-|\tilde{\pi}_{m}|)}.
\end{equation}  
In the case of $\pi_{m}(\alpha_{1}*\cdots*\alpha_{m})$, the situation is 
almost the same. The only difference is that the number of added parentheses 
increases by $|\pi_{m}|-|\tilde{\pi}_{m}|$. Therefore, by assumption of 
induction, the coefficient of $\pi_{m}(\alpha_{1}*\cdots*\alpha_{m})$
equals $(-1)^{(m-|\tilde{\pi}_{m}|)}\cdot
(-1)^{-(|\pi_{m}|-|\tilde{\pi}_{m}|)}=(-1)^{(m-|\pi_{m}|)}$. Thus, (\ref{ql})
is derived for $l=m$, and the proof of Proposition 1 is completed. Q.E.D.   

And the proof of the Theorem 4 is completed.                  Q.E.D.
\\
\\
{\bf Proof of Theorem 1 )}\\
From the statement of Theorem 4, we obtain the formulas in the 
case of $(k+2\leq N \leq 2k)$:
\begin{eqnarray}
&&\gamma_{0}^{N,k,1}(w)=\gamma_{0}^{2k,k,1}(w)=k\prod_{j=1}^{k-1}(k+jw),\no\\
&&\gamma_{0}^{N,k,d}(w)=(\prod_{j=1}^{d-1}(1+jw))^{2k-N}
\gamma_{0}^{2k,k,d}(w)=0,\;\;\;\;(d\geq 2).
\label{gam}
\end{eqnarray}
But from the proof of the Corollary 1, we can determine 
$L^{N,k,d}_{n}$ only using (\ref{gam}). Therefore, we can conclude that 
the recursive formulas in Theorem 1 compute $L^{N,k,d}_{n}$ 
correctly.   Q.E.D. \\
\subsection{$N-k=1$ case} 
In this case, we had better introduce $\tilde{\psi}_{\alpha}(t)$,
\begin{equation}
\tilde{\psi}_{\alpha}(t):=\exp(k!\cdot\exp(t))\cdot\psi_{\alpha}(t),\;\;\;
(\alpha=0,1,\cdots,N-2),
\label{B}
\end{equation}
instead of $\psi_{\alpha}(t)$ in (\ref{gm1})
because of the following Theorem of Givental \cite{giv}.
\begin{theorem}
{\bf (Givental)}\\
If $N-k=1$, $\tilde{\psi}_{0}(t)$  satisfy the rank $N-1$ ODE:
\begin{equation}
\biggl((\partial_{t})^{N-1}-k\cdot e^{t}\cdot(k\partial_{t}+k-1)\cdots(k\partial_{t}+2)\cdot(k\partial_{t}+1)\biggr)w(t)=0.
\label{simB} 
\end{equation}
\end{theorem}
This theorem is equivalent to $L_{m}^{k+1,k,1}=\tilde{L}_{m}^{k+1,k,1}-k!,\;\;
L_{m}^{k+1,k,d}=\tilde{L}_{m}^{k+1,k,d}\;\;(d\geq 2)$ \cite{cj}.  
\subsection{Calabi-Yau case ($N-k=0$)}
Then we turn into the case of Calabi-Yau hypersurface. To clarify the meaning 
of the virtual structure constants introduced in
\cite{cj}, we had better introduce the B-model deformation parameter $x$
instead of $t$ and consider the following Gauss-Manin system. 
\begin{eqnarray}
&&\d_{x}\tilde{\psi}_{-1}(x)=\tilde{L}^{k,k}_{k-1}(e^{x})\cdot\tilde{\psi}_{0}(x),\no\\
&&\d_{x}\tilde{\psi}_{n}(x)=\tilde{L}^{k,k}_{k-2-n}(e^{x})\cdot\tilde{\psi}_{n+1}(x),\;\;\;(n=0,\cdots..k-3)
\no\\
&&\d_{x}\tilde{\psi}_{k-2}(x)=\tilde{L}^{k,k}_{0}(e^{x})\cdot\tilde{\psi}_{k-1}(x).
\label{vg0}
\end{eqnarray}
We can derive the following equality from the above equations:
\begin{equation}
\tilde{\psi}_{k-1}(x)=\frac{1}{\tilde{L}_{0}^{k,k}(e^x)}(\partial_{x}(\frac{1}{\tilde{L}_{1}^{k,k}(e^x)}\cdots \partial_{x}(\frac{1}{\tilde{L}_{k-2}^{k,k}(e^x)}\partial_{x}(\frac{1}{\tilde{L}_{k-1}^{k,k}(e^x)}\d_{x}\tilde{\psi}_{-1}(x)))\cdots)).
\label{virtual}
\end{equation}
(\ref{virtual}) motivates us to state the following theorem. 
\begin{theorem}
\begin{eqnarray}
&&\frac{1}{\tilde{L}_{0}^{k,k}(e^x)}(\partial_{x}(\frac{1}{\tilde{L}_{1}^{k,k}(e^x)}\cdots \partial_{x}(\frac{1}{\tilde{L}_{k-2}^{k,k}(e^x)}\partial_{x}(\frac{1}{\tilde{L}_{k-1}^{k,k}(e^x)}w(x)))\cdots))\no\\
&&=\biggl((\partial_{x})^{k-1}-k\cdot e^{x}\cdot(k\partial_{x}+k-1)\cdots(k\partial_{x}+2)\cdot(k\partial_{x}+1)\biggr)w(x)
\label{pfN}
\end{eqnarray}
\end{theorem}
{\it proof)} We only have to apply formally the discussion of $N-k\geq 2$ 
case to the $N-k=0$ case with the Gauss-Manin system (\ref{vg0}). Q.E.D.

Since $\tilde{L}_{n}^{k,k}(e^x)=\tilde{L}_{k-1-n}^{k,k}(e^x)$, we have,
\begin{cor}
\begin{equation}
u_{j}^{k,k}(x):={\tilde{L}_{0}^{k,k}(e^x)}\int^{x}dx_{1}
{\tilde{L}_{1}^{k,k}(e^{x_{1}})}\int^{x_{1}}dx_{2}
{\tilde{L}_{2}^{k,k}(e^{x_{2}})}\cdots\int^{x_{j-1}}dx_{j}
{\tilde{L}_{j}^{k,k}(e^{x_{j}})}.
\label{cy}
\end{equation}
\end{cor}
\begin{Rem}
Representation of the Picard-Fuchs differential equation given in (\ref{pfN})
can also be seen in \cite{bara}. We think that our approach via Gauss-Manin 
system is a kind of reduction of the method used in \cite{bara}
, restricted to the K\"ahler deformation. 
\end{Rem}
(\ref{cy}) enables us to write out $\tilde{L}_{n}^{k,k}(e^{x})$ explicitly 
in terms of the solution of the Picards-Fuchs differentail 
equation used in the Mirror computation in \cite{gmp}, \cite{mnmj}. For example, we have:
\begin{eqnarray}
\tilde{L}_{0}^{k,k}(e^{x})
&=& w_{0}^{k,k}(x),\\
\tilde{L}_{1}^{k,k}(e^{x})&=&\d_{x}(x+\frac{w_{1}^{k,k}(x)}
{w_{0}^{k,k}(x)}),\\
\tilde{L}_{2}^{k,k}(e^{x})&=&\d_{x}
(x+\frac{2w_{1}^{k,k}(x)w_{0}^{k,k}(x)+\d_{x}w_{2}^{k,k}(x)w_{0}^{k,k}(x)-
w_{2}^{k,k}(x)\d_{x}w_{0}^{k,k}(x)}
{2((w_{0}^{k,k}(x))^{2}+\d_{x}w_{1}^{k,k}(x)w_{0}^{k,k}(x)-
w_{1}^{k,k}(x)\d_{x}w_{0}^{k,k}(x))}).
\end{eqnarray}
These results agree with the computation in \cite{gmp}, and they give 
us the proof of Theorem 2.
\subsection{Extension to the General Type Hypersurfaces ($N-k<0$)}
If $N-k<0$, We consider the rank $k-1$ ODE:
\begin{eqnarray}
&&\biggl((\partial_{x})^{N-1}-k\cdot e^{x}\cdot(k\partial_{x}+k-1)\cdots(k\partial_{x}+2)\cdot(k\partial_{x}+1)\biggr)w(x)\no\\
&&\biggl(1-k\cdot e^{x}\cdot(k\partial_{x}+k-1)\cdots(k\partial_{x}+2)\cdot(k\partial_{x}+1)\frac{1}{(\d_{x})^{N-1}}\biggr)\d_{x}^{N-1}w(x)\no\\
&&=0.
\label{neg} 
\end{eqnarray}
Here, we propose the B-model Gauss-Manin system associated to (\ref{neg}):
\begin{equation}
\partial_{x}\tilde{\psi}_{\alpha}(x)=\tilde{\psi}_{\alpha+1}(x)+
\sum_{d=1}^{\infty}\exp(dx)\cdot\tilde{L}^{N,k,d}_{N-2-\alpha}\cdot\tilde
{\psi}_{\alpha+1+(k-N)d}(x),
\label{gen}
\end{equation}
where $\alpha$ runs through ${\bf Z}$. $\tilde{L}^{N,k,d}_{n}$ is the virtual 
structure constant introduced in \cite{gene}.
$\tilde{L}^{N,k,d}_{n}$ is non-zero 
if $0\leq n \leq N-1+(k-N)d$ and if $d\geq 1$. All the non-vanishing $\tilde{L}^{N,k,d}_{n}$'s are evaluated via the recursive formulas proposed in 
\cite{6}.  
 Therefore, we have infinite number of non-vanishing virtual structure 
constants in this case. Straightforward application of the discussion of 
$N-k\geq 2$ case to this case leads us to the following theorem:
\begin{theorem}
We can reconstruct the ODE (\ref{neg}) from (\ref{gen}). Conversely, we can 
determine $\tilde{L}^{N,k,d}_{n}$ by (\ref{neg}). 
\end{theorem}
Now, we explain the reconstruction process of (\ref{neg}) from (\ref{gen}).
First we introduce the algebra of differential operator $\d_{x}$.
\begin{eqnarray}
&&\d_{x}\cdot\frac{1}{\d_{x}}=\frac{1}{\d_{x}}\cdot\d_{x}=1,\no\\
&&\d_{x}e^{jx}=e^{jx}(\d_{x}+j),\;\;\;\;
\frac{1}{\d_{x}}e^{jx}=e^{jx}\frac{1}{(\d_{x}+j)}.
\label{diff}
\end{eqnarray}
Using (\ref{diff}), we can obtain the following formula:
\begin{eqnarray}
&&\biggl(1-k\cdot e^{x}\cdot(k\partial_{x}+k-1)\cdots(k\partial_{x}+2)\cdot
(k\partial_{x}+1)\frac{1}{(\d_{x})^{N-1}}\biggr)^{-1}\no\\
&&=1+\sum_{d=1}^{\infty}e^{dx}\cdot\prod_{m=0}^{kd-1}
(k\d_{x}+m)\prod_{j=0}^{d-1}\frac{1}{(\d_{x}+j)^{N}}.
\end{eqnarray}
Looking back at (\ref{gen}), we can easily see,
\begin{equation}
\tilde{\psi}_{j}(x)=(\d_{x})^{j-N+1}\tilde{\psi}_{N-1}(x),\;\;\; (j\geq N-1).
\end{equation}  
Then using the algebras in (\ref{diff}) and (\ref{gen}), we can 
inductively construct the 
pseudo-differential operator $F_{j}(e^{x},\d_{x})$ when $j\geq 0 $,
\begin{equation}
(\d_{x})^{N-1+j}\tilde{\psi}_{-j}(x)=F_{j}(e^{x},\d_{x})\tilde{\psi}_{N-1}(x).
\label{pd}
\end{equation}
Then we consider the limit $F_{\infty}(e^{x},\d_{x}):=\lim_{j\rightarrow \infty}F_{j}(e^{x},\d_{x})$. Now, our assertion is the following statement:
\begin{equation}
\frac{1}{F_{\infty}(e^{x},\d_{x})}(\d_{x})^{N-1}=(1-k\cdot e^{x}\cdot(k\partial_{x}+k-1)\cdots(k\partial_{x}+2)\cdot(k\partial_{x}+1)\frac{1}{(\d_{x})^{N-1}})
(\d_{x})^{N-1},
\end{equation}
or equivalently, 
\begin{equation}
F_{\infty}(e^{x},\d_{x})=1+\sum_{d=1}^{\infty}e^{dx}\cdot\prod_{m=0}^{kd-1}
(k\d_{x}+m)\prod_{j=0}^{d-1}\frac{1}{(\d_{x}+j)^{N}}.
\end{equation}
Conversely, we can determine $\tilde{L}_{n}^{N,k,d}$ assuming the above 
equation. This process corresponds to the 
B-model computation in the $N=k$ case,
and combining it with generalized mirror transformation, we can construct 
``mirror computation'' to the general type hypersurface $M_{N}^{k}\;\;(N<k)$.
\\     
{\bf Acknowlegement}\\
The author especially thanks Prof. B. Kim for valuable discussions and 
for invitation to Korea 
Institute for Advanced Study, where part of this work was done.
He also thanks the organizers of ``Summer Institute 2001'' and of 
``International Workshop on Integrable Models, Combinatorics and 
Representation Theory'' for the hospitality during the finishing 
period of this work.
Research of the author is partially supported by the grant of Japan 
Society for Promotion of Science.

\end{document}